\documentclass{amsart}
\usepackage{latexsym}
\usepackage{amssymb}
\usepackage{amsmath}
\usepackage{amsfonts}
\usepackage{amsthm}
\usepackage{amscd}
\usepackage{texdraw}
\usepackage{color}
\usepackage{enumitem}
\usepackage[centertags]{amsmath}

\newlength{\defbaselineskip}
\newcommand{\setlinespacing}[1]%
           {\setlength{\baselineskip}{#1 \defbaselineskip}}

\newcommand{\vertiii}[1]{{\left\vert\kern-0.25ex\left\vert\kern-0.25ex\left\vert #1
    \right\vert\kern-0.25ex\right\vert\kern-0.25ex\right\vert}}
\theoremstyle{definition}
\newtheorem{thm}{Theorem}
\newtheorem{prop}{Proposition}
\newtheorem{cor}{Corollary}
\newtheorem{lem}{Lemma}

\newtheorem{remark}{Remark}

\allowdisplaybreaks
\begin{document}
\title[Analogues of Wiener's Theorem]{Vector valued Beurling algebra analogues of Wiener's Theorem}

\author[P. A. Dabhi]{Prakash A. Dabhi}
\address{Institute of Infrastructure Technology Research and Management(IITRAM), Maninagar (East), Ahmedabad - 380026, Gujarat, India}
\email{lightatinfinite@gmail.com, prakashdabhi@iitram.ac.in}

\author[K. B. Solanki]{Karishman B. Solanki}
\address{Institute of Infrastructure Technology Research and Management(IITRAM), Maninagar (East), Ahmedabad - 380026, Gujarat, India}
\email{karishsolanki002@gmail.com}

\date{}

\begin{abstract}
Let $0<p\leq 1$, $\omega$ be a weight on $\mathbb Z$, and let $\mathcal A$ be a unital Banach algebra. If $f$ is a continuous function from the unit circle $\mathbb T$ to $\mathcal A$ such that $\sum_{n\in \mathbb Z} \|\widehat f(n)\|^p \omega(n)^p<\infty$ and $f(z)$ is left invertible for all $z \in \mathbb T$, then there is a weight $\nu$ on $\mathbb Z$ and a continuous function $g:\mathbb T \to \mathcal A$ such that $1\leq \nu \leq \omega$, $\nu$ is constant if and only if $\omega$ is constant, $g$ is a left inverse of $f$ and $\sum_{n\in \mathbb Z}\|\widehat g(n)\|^p\nu(n)^p<\infty$. We shall obtain a similar result when $\omega$ is an almost monotone algebra weight and $1<p<\infty$. We shall obtain an analogue of this result on the real line. We shall apply these results to obtain $p-$power weighted analogues of the results of off diagonal decay of infinite matrices of operators.
\end{abstract}



\maketitle

\section{Introduction}
We shall use the following notations: $\mathbb{N}$ is the set of all positive integers, $\mathbb{N}_0=\mathbb{N}\cup\{0\}$, $\mathbb{Z}$ is the set of all integers, $\mathbb{R}$ is the set of all real numbers, $\mathbb{C}$ is the set of all complex numbers, $\mathbb{T}=\{z\in\mathbb{C}:|z|=1\}$ is the unit circle in $\mathbb{C}$, and for $0\leq a\leq b$, $\Gamma(a,b)=\{z\in\mathbb{C}:a\leq |z|\leq b\}$.

A map $\omega:\mathbb Z \to [1,\infty)$ is a \emph{weight} if $\omega(m+n)\leq \omega(m)\omega(n)$ for all $m,n \in \mathbb Z$. For a weight $\omega$ on $\mathbb{Z}$, let $\rho_{1,\omega}=\sup\{\omega(-n)^{-1/n}:n \in \mathbb{N}\}$ and $\rho_{2,\omega}=\inf\{\omega(n)^{1/n}:n \in \mathbb{N}\}$. Then $\rho_{1,\omega} \leq 1 \leq \rho_{2,\omega}$. A weight $\omega$ on $\mathbb{Z}$ is \emph{admissible} if $\rho_{1,\omega}=1=\rho_{2,\omega}$. Let $1<p<\infty$, and let $q$ be its conjugate index, i.e., $\frac{1}{p}+\frac{1}{q}=1$. A weight $\omega$ on $\mathbb{Z}$ is an \emph{almost monotone algebra weight} if the following conditions hold:
\begin{enumerate}
\item $\omega^{-q}\star\omega^{-q}\leq\omega^{-q}$ and $\sum_{n \in \mathbb{Z}} \omega(n)^{-q}<\infty$.
\item If $\rho_{1,\omega}=1$, then there is a positive constant $K$ such that $\omega(n)\leq K\omega(n+k)$ for all $-n,-k \in \mathbb{N}$.
\item If $\rho_{2,\omega}=1$, then there is a positive constant $K$ such that $\omega(n)\leq K\omega(n+k)$ for all $n,k \in \mathbb{N}_0$.
\end{enumerate}

Let $f:\mathbb{T}\to\mathbb{C}$ be a nowhere vanishing continuous function. The series $\sum_{n\in\mathbb{Z}}\widehat{f}(n)e^{int}$ is the \emph{Fourier series} of $f$, where $\widehat{f}(n)=\frac{1}{2\pi}\int_{0}^{2\pi}f(e^{it})e^{-int}dt \; (n\in\mathbb{Z})$ are the \emph{Fourier coefficients} of $f$.

If $f$ has absolutely convergent Fourier series, i.e., $\sum_{n\in\mathbb{Z}}|\widehat{f}(n)|<\infty$, then the classical Wiener's theorem \cite{Wi} states that the reciprocal $1/f$ also has absolutely convergent Fourier series.

Let $\omega$ be a weight on $\mathbb{Z}$, and let $p>0$. The function $f$ has \emph{$\omega$-absolutely convergent Fourier series ($\omega$-ACFS)} if $\sum_{n\in\mathbb{Z}}|\widehat f(n)|\omega(n)<\infty$ and $f$ has \emph{$p$-th power $\omega$-absolutely convergent Fourier series ($p\omega$-ACFS)} if $\sum_{n\in\mathbb{Z}}|\widehat f(n)|^p\omega(n)^p<\infty$.

For $0<p\leq1$, $p$-power analogue of Wiener's theorem was obtained by \.Zelazko \cite{Ze}. Let $\omega$ be a non-quasianalytic weight on $\mathbb{Z}$, i.e., $\sum_{n\in\mathbb{Z}}\frac{\log\omega(n)}{1+n^2}<\infty$. If $f$ has $\omega$-ACFS, then Domar \cite{Do} proved that $1/f$ also has $\omega$-ACFS. Note that non-quasianalytic weights are admissible. Let $\omega$ be any weight on $\mathbb{Z}$. If $f$ has $\omega$-ACFS, then Bhatt and Dedania \cite{Bh} constructed a weight $\nu$ on $\mathbb{Z}$ such that $1\leq\nu\leq\omega,\ \nu$ is constant if and only if $\omega$ is constant and $1/f$ has $\nu$-ACFS. An analogue of Wiener's theorem in weighted case for $0<p\leq 1$ is obtained in \cite{De}. Let $1<p<\infty$, and let $\omega$ be an almost monotone algebra weight on $\mathbb{Z}$. If $f$ has $p\omega$-ACFS, then an existence of an almost monotone algebra weight $\nu$ on $\mathbb{Z}$ such that $1\leq\nu\leq K\omega$ for some $K>0$ and $1/f$ has $\nu$-ACFS is obtained in \cite{Da}.

Note that all these results deal with the complex valued functions on the unit circle $\mathbb{T}$ and the condition $f$ is nowhere vanishing implies that the complex number $f(z)$ is invertible in $\mathbb{C}$ for all $z\in\mathbb{T}$. Let $\mathcal A$ be a unital Banach algebra, and let $f:\mathbb T \to \mathcal A$ be a continuous function. The series $\sum_{n\in \mathbb Z} \widehat f(n) e^{int}$ is the \emph{Fourier series} of $f$ and $\widehat f(n)=\frac{1}{2\pi}\int_0^{2\pi}f(e^{it})e^{-int}dt \in \mathcal A\;(n\in \mathbb Z)$ are the \emph{Fourier coefficients} of $f$. If $f$ has a representation $f(z)=\sum_{n\in \mathbb Z}\widehat f(n)z^n \;(z\in \mathbb T)$ such that $\sum_{n\in \mathbb Z}\|\widehat f(n)\|<\infty$ and if $f(z)$ is left invertible (respectively, right invertible, invertible) for all $z \in \mathbb T$, then Bochner and Phillips in \cite{Bo} proved that a left inverse (respectively, right inverse, inverse) $g$ of $f$ has a representation $g(z)=\sum_{n\in \mathbb Z} \widehat g(n) z^n \;(z\in \mathbb T)$ and $\sum_{n\in \mathbb Z}\|\widehat g(n)\|<\infty$. Let $\omega$ be an admissible weight on $\mathbb{Z}$. In \cite{Ba}, Baskakov proved that if $f$ has $\omega$-ACFS, i.e., $\sum_{n\in \mathbb Z}\|\widehat f(n)\|\omega(n)<\infty$ and if $f(z)$ is invertible in $\mathcal{A}$ for all $z\in\mathbb{T}$, then the inverse of $f$ exists, say $g$, and $g$ also has $\omega$-ACFS. We obtain weighted analogues of these results without admissibility of weight. All algebras considered in this paper are complex algebras.

Let $\mathcal{A}$ be a Banach algebra. The norm of $\mathcal A$ will be denoted by $\|\cdot\|$ throughout the paper except where stated otherwise. By $C(\mathbb T,\mathcal A)$, denote the collection of all continuous functions from $\mathbb T$ to $\mathcal A$. Let $0<p<\infty$, and let $\omega$ be a weight on $\mathbb{Z}$. A function $f \in C(\mathbb T,\mathcal A)$ has \emph{$p$-th power $\omega$-absolutely convergent Fourier series ($p\omega$- ACFS)} if $f(z)=\sum_{n\in \mathbb Z}\widehat f(n)z^n\;(z \in \mathbb T)$ and $\sum_{n \in \mathbb Z}\|\widehat f(n)\|^p\omega(n)^p<\infty$.

We shall prove the following results.
\begin{thm} \label{thm1}
Let $0<p\leq 1$. Let $\mathcal{A}$ be a unital Banach algebra, $\omega$ be a weight on $\mathbb{Z}$, $f \in C(\mathbb{T},\mathcal{A})$, and let $f$ have $p\omega$-ACFS. If $f(z)$ is left invertible (respectively, right invertible, invertible) for all $z \in \mathbb{T}$, then there exist a weight $\nu$ on $\mathbb{Z}$ and $g \in C(\mathbb{T},\mathcal{A})$ such that
\begin{enumerate}
\item $1 \leq \nu \leq \omega$;
\item $\nu$ is constant if and only if $\omega$ is constant;
\item $g$ has $p\nu$-ACFS;
\item $g$ is a left inverse (respectively, right inverse, inverse) of $f$.
\end{enumerate}
In particular, if $\omega$ is an admissible weight, then a left inverse (respectively, right inverse, inverse) of $f$ has $p\omega$-ACFS.
\end{thm}

\begin{thm} \label{thm2}
Let $1<p<\infty$, $q$ be the conjugate index of $p$, and let $\omega$ be an almost monotone algebra weight on $\mathbb Z$. Let $\mathcal{A}$ be a unital Banach algebra, $f \in C(\mathbb{T},\mathcal{A})$, and let $f$ have $p\omega$-ACFS. If $f(z)$ is left invertible (respectively, right invertible, invertible) for all $z \in \mathbb{T}$, then there exist an almost monotone algebra weight $\nu$ on $\mathbb{Z}$ and $g \in C(\mathbb{T},\mathcal{A})$ such that
\begin{enumerate}
\item $1 \leq \nu \leq K\omega$ for some $K>0$;
\item $g$ has $p\nu$-ACFS;
\item $g$ is a left inverse (respectively, right inverse, inverse) of $f$.
\end{enumerate}
In particular, if $\omega$ is an admissible weight, then a left inverse (respectively, right inverse, inverse) of $f$ has $p\omega$-ACFS.
\end{thm}

A Borel measurable map $\omega:\mathbb R \to [1,\infty)$ is a \emph{weight} if $\omega(x+y)\leq \omega(x)\omega(y)$ for all $x,y \in \mathbb R$. For a weight $\omega$ on $\mathbb{R}$, let $\rho_{1,\omega}=\sup\{\log\omega(x)^{1/x}:x<0\}$ and $\rho_{2,\omega}=\inf\{\log\omega(x)^{1/x}:x>0\}$. Then $\rho_{1,\omega} \leq 0 \leq \rho_{2,\omega}$. A weight $\omega$ on $\mathbb{R}$ is \emph{admissible} if $\rho_{1,\omega}=0=\rho_{2,\omega}$. Let $1<p<\infty$, and let $q$ be the conjugate index of $p$. A weight $\omega$ on $\mathbb{R}$ is an \emph{almost monotone algebra weight} if it satisfies the following conditions:
\begin{enumerate}
\item $\omega^{-q}\star\omega^{-q}\leq\omega^{-q}$ and $\int_{\mathbb{R}}\omega^{-q}<\infty$.
\item If $\rho_{1,\omega}=0$, then there is a positive constant $K$ such that $\omega(x)\leq K\omega(x+y)$ for all $x,y\leq0$.
\item If $\rho_{2,\omega}=0$, then there is a positive constant $K$ such that $\omega(x)\leq K\omega(x+y)$ for all $x,y\geq0$.
\end{enumerate}

Let $\mathcal{A}$ be a unital Banach algebra with unit $\mathbf{e}$, $\omega$ be a weight on $\mathbb{R}$, $1\leq p<\infty$, and let  $$L^p(\mathbb{R},\omega,\mathcal{A})=\{f:\mathbb{R}\to\mathcal{A}:\int_\mathbb{R}\|f(x)\|^p\omega(x)^pdx<\infty\}.$$ By $L^p(\mathbb{R},\omega,\mathcal{A})_\mathbf1$, we denote the unitization of $L^p(\mathbb{R},\omega,\mathcal{A})$ by adjoining the unit $\mathbf1$ when $L^p(\mathbb{R},\omega,\mathcal{A})$ is an algebra under convolution.

In \cite{Bo}, Bochner and Phillips proved that if $f\in L^1(\mathbb{R},\mathcal{A})$ and if $\mathbf{e}+\int_\mathbb{R}f(x)e^{ixt}dt$ has a left inverse for all $t\in\mathbb{R}$, then there is $g\in L^1(\mathbb{R},\mathcal{A})$ such that $\mathbf1+g$ is a left inverse of $\mathbf1+f$. The following theorems are $p$-power weighted analogues of this result and they may be seen as continuous analogue of the above results.

\begin{thm}\label{thr}
Let $\mathcal{A}$ be a unital Banach algebra with unit $\mathbf e$, $\omega$ be a weight on $\mathbb{R}$, and let $f\in L^1(\mathbb{R},\omega,\mathcal{A})$ be such that $\mathbf e+\int_\mathbb{R} f(x) e^{itx}dx$ is left invertible (respectively, right invertible, invertible) in $\mathcal{A}$ for all $t\in\mathbb{R}$. Then there is a weight $\nu$ on $\mathbb{R}$ and $g\in L^1(\mathbb{R},\nu,\mathcal{A})$ such that $1\leq\nu\leq\omega$, $\nu$ is constant if and only if $\omega$ is constant and $\mathbf1+g$ is a left inverse (respectively, right inverse, inverse) of $\mathbf1+f$. In particular, if $\omega$ is admissible, then $\mathbf{1}+f$ is left invertible (respectively, right invertible, invertible) in $L^1(\mathbb{R},\omega,\mathcal{A})_\mathbf1$.
\end{thm}

\begin{thm}\label{thrp}
Let $1<p<\infty$, and let $q$ be its conjugate index. Let $\mathcal{A}$ be a unital Banach algebra with unit element $\mathbf e$, $\omega$ be an almost monotone algebra weight on $\mathbb{R}$, and let $f\in L^p(\mathbb{R},\omega,\mathcal{A})$ be such that $\mathbf e+\int_\mathbb{R} f(x) e^{itx}dx$ is left invertible (respectively, right invertible, invertible) in $\mathcal{A}$ for all $t\in\mathbb{R}$. Then there is an almost monotone algebra weight $\nu$ on $\mathbb{R}$ and $g\in L^p(\mathbb{R},\nu,\mathcal{A})$ such that $1\leq\nu\leq K \omega$ for some $K>0$ and $\mathbf1+g$ is a left inverse (respectively, right inverse, inverse) of $\mathbf1+f$. In particular, if $\omega$ is admissible, then $\mathbf{1}+f$ is left invertible (respectively, right invertible, invertible) in $L^p(\mathbb{R},\omega,\mathcal{A})_\mathbf1$.
\end{thm}

Let $\mathcal A$ and $\mathcal B$ be Banach algebras, $\mathcal A \subset \mathcal B$, and let $\mathcal A$ and $\mathcal B$ have the same unit. Then $\mathcal A$ is \emph{inverse closed} (\emph{spectrally invariant}) in $\mathcal B$ if $a \in \mathcal A$ and $a^{-1}\in \mathcal B$ imply $a^{-1}\in \mathcal A$. Wiener's theorem states that the Banach algebra of continuous functions on $\mathbb T$ having absolutely convergent Fourier series is inverse closed in the Banach algebra of all continuous functions on $\mathbb T$. For various generalization of Wiener's Theorem we refer \cite{Ra, Kr, Bo, Fa, Ge, GrMc, Gro1, GrLe, gr, gr1, Kr1, Na, Sh, Su, Ta} and references therein, of course the list is not exhaustive. Let $X$ be a complex Banach space and $B(X)$ be the Banach algebra of all bounded linear operators on $X$. A map $R: \mathbb{Z} \to B(X)$ is a \emph{resolution of the identity} if $R(n)$ is a projection on $X$ for all $n$, $R(k)R(l)=0$ for all $k \neq l$ and the series $\sum_{k \in \mathbb{Z}} R(k)x$ converges to $x$ unconditionally for all $x \in X$. The resolution of identity satisfies some additional conditions (defined in section 6). For each $A \in B(X)$, consider the matrix $(A_{ij})_{i,j \in \mathbb{Z}}$, where $A_{ij}=R(i)AR(j) \in B(X)$. Define $d_A: \mathbb{Z} \to [0,\infty)$ by $d_A(k)=\sup_{i-j=k}\|A_{ij}\| \; (k \in \mathbb{Z})$. Then $d_A$ characterizes the off-diagonal decay of the entries of $(A_{ij})$. Let $\omega$ be an admissible weight on $\mathbb Z$. Let $\mathcal B$ be the collection of all $A \in B(X)$ such that $\|A\|=\sum_{n\in \mathbb Z}d_A(n)\omega(n)<\infty$. Then $(\mathcal B,\|\cdot\|)$ is a unital Banach algebra having unit same as that of $B(X)$. Baskakov proved in \cite{Ba1} that $\mathcal B$ is inverse closed in $B(X)$. Let $\Lambda$ be a lattice in $\mathbb R^d$, and let $\omega$ be an admissible weight on $\Lambda$. For a matrix $A=(a_{\lambda,\mu})_{\lambda,\mu \in \Lambda}$, let $d_A(\mu)=\sup_{\lambda \in \Lambda}|a_{\lambda,\lambda-\mu}|$. Let $\mathcal B$ be the collection of all $A=(a_{\lambda,\mu})$ such that $A$ defines a bounded linear operator on $\ell^2(\Lambda)$ and $\|A\|=\sum_{\mu\in \Lambda}d_A(\mu)\omega(\mu)<\infty$. Gr\"ochenig and Rzeszotnik in \cite{G1R} proved that $\mathcal B$ is inverse closed in $B(\ell^2(\Lambda))$ and obtained some more general results. 

Let $\omega$ be a weight on $\mathbb Z$, not necessarily admissible. If $A$ is in any of the Banach algebras mentioned above and if $A$ is invertible in $B(\mathcal H)$ or $B(\ell^2(\mathbb Z))$, then the matrix of $A^{-1}$ may not satisfy the same weighted off-diagonal decay. Let $0<p<\infty$. Let $\omega$ be a weight on $\mathbb Z$. When $p>1$, we shall require $\omega$ to be an almost monotone algebra weight. Let $R$ be a resolution of the identity on the Banach space $X$ with certain conditions. Let $\mathcal A_{p,\omega}$ be the collection of all $A=(A_{i,j}) \in B(X)$ such that $\sum_{n\in\mathbb{Z}}d_A(n)^p\omega(n)^p<\infty$. Then $\mathcal A_{p,\omega}$ will be a Banach algebra ($p$-Banach algebra when $0<p\leq1$) with appropriate norm and is contained in $B(X)$ having the same identity. The algebra $\mathcal A_{p,\omega}$ is not necessarily inverse closed in $B(X)$. If $A$ is invertible in $B(X)$, then the sum $\sum_{n\in \mathbb Z}d_{A^{-1}}(n)^p\omega(n)^p$ is not necessarily finite. When $A \in \mathcal A_{p,\omega}$ is invertible in $B(X)$, we find a weight $\nu$ which is closely related to $\omega$ such that $A$ has inverse in $\mathcal A_{p,\nu}$. In the case when $\omega$ is admissible, we get $\nu=\omega$ and hence $A^{-1}\in \mathcal A_{p,\omega}$. Our results generalize the results in \cite{Ba} and \cite{Bo} in weighted cases and for $0<p<\infty$.


\section{Proof of Theorem \ref{thm1}}
Let $0<p\leq1$, and let $\mathcal{A}$ be an algebra. A mapping $\|\cdot\| : \mathcal{A} \to [0,\infty)$ is a \emph{$p$-norm} \cite{De} on $\mathcal{A}$ if the following conditions hold for all $x,y\in \mathcal{A}$ and $\alpha \in \mathbb{C}$.
\begin{enumerate}
\item $\|x\|=0$ if and only if $x=0$;
\item $\|x+y\|\leq\|x\|+\|y\|$;
\item $\|\alpha x\|= |\alpha|^p \|x\|$;
\item $\|xy\|\leq\|x\|\|y\|$.
\end{enumerate}
If $\mathcal{A}$ is complete in the $p$-norm, then $(\mathcal{A},\|\cdot\|)$ is a \emph{$p$-Banach algebra} \cite{De}. When $p=1$, $\mathcal A$ is a Banach algebra.

Let $\mathcal{A}$ be an algebra without unit element, and let $\mathcal{A}_\mathbf e=\mathcal{A}\times\mathbb{C}=\{(a,\alpha):a\in\mathcal{A}, \alpha\in\mathbb{C}\}$. Then $\mathcal{A}_\mathbf e$ is an algebra with pointwise linear operations  and the multiplication $(a,\alpha)(b,\beta)=(ab+\alpha b+\beta a,\alpha\beta)$ for all $a,b\in\mathcal{A}$ and $\alpha,\beta\in\mathbb{C}$. The element $\mathbf e=(0,1)$ is the unit element of $\mathcal{A}_\mathbf e$. If $\mathcal{A}$ is a $p$-Banach algebra, then $\mathcal{A}_\mathbf e$ is also a $p$-Banach algebra with the norm $$\|(a,\alpha)\|=\|a\|+|\alpha|^p \quad (a\in\mathcal{A},\alpha\in\mathbb{C}).$$ Since $\mathcal{A}$ can be identified with the ideal $\{(a,0):a\in\mathcal{A}\}$, of codimension 1, of $\mathcal{A}_\mathbf e$, the element $(a,\alpha)$ of $\mathcal{A}_\mathbf e$ is customarily written as $a+\alpha \mathbf e$ or $\alpha \mathbf e+a$. The algebra $\mathcal{A}_\mathbf e$ is the \emph{unitisation} of $\mathcal{A}$.

Let $\mathcal{A}$ be a commutative $p$-Banach algebra. A nonzero linear map $\varphi:\mathcal A \to \mathbb C$ satisfying $\varphi(ab)=\varphi(a)\varphi(b)\;(a,b \in \mathcal A)$ is a \emph{complex homomorphism} on $\mathcal A$. Let $\Delta(\mathcal A)$ be the collection of all complex homomorphisms on $\mathcal A$. For $a \in \mathcal A$, let $\widehat a:\Delta(\mathcal A)\to \mathbb C$ be $\widehat a(\varphi)=\varphi(a)\;(\varphi\in \Delta(\mathcal A))$. The smallest topology on $\Delta(\mathcal A)$ making each $\widehat a$, $a\in \mathcal A$, continuous is the \emph{Gel'fand topology} on $\Delta(\mathcal A)$ and $\Delta(\mathcal A)$ with the Gel'fand topology is the \emph{Gel'fand space} of $\mathcal{A}$. For more details on it refer \cite{Ge, Ze}.

Let $X$ be a complex Banach space, and let $B(X)$ be the collection of all bounded linear transformation from $X$ to itself. If $X_0\subset X$ and $A_0\subset B(X)$, then $A_0X_0=\{T(x)\in X:T\in A_0, x\in X_0\}$. We say that $X$ is \emph{irreducible} over $A_0$ if for every non-zero $x\in X, \ A_0x=X$.

\begin{lem} \cite[Lemma 9]{Bo} \label{lbo9}
If $X$ is irreducible over $A_0$, and if $A'=\{T\in B(X) : ST=TS \ \text{for all} \ S\in A_0\}$, then $A'$ is isomorphic to $\mathbb{C}$.
\end{lem}

\begin{lem}\cite[Lemma 10 and 11]{Bo}\label{lbo10}
If $\mathcal{B}$ is a unital algebra, $I$ is a maximal left ideal, $V$ is addition group of the cosets of $\mathcal{B}/I$, and if $B_0=\{L_x:x\in\mathcal{B}\}$, where $L_x:V\to V$ are defined as $L_x(v)=xv \quad (v\in V)$, then $V$ is irreducible over $B_0$. Also, if for a fixed $x\in\mathcal{B}$ and every maximal left ideal $I$, the corresponding element $L_x$ of $B_0$ has a left inverse in $B_0$, then $x$ has a left inverse in $\mathcal{B}$.
\end{lem}

\begin{lem}\label{lbo13}
If $\mathcal{B},I,V$ and $B_0$ are as in Lemma \ref{lbo10} with $\mathcal{B}$ being $p$-Banach algebra, then $V$ is a $p$-Banach space, $B_0$ is a $p$-normed algebra and $\|L_x\| \leq \|x\| \quad (x\in\mathcal{B})$, where $\|L_x\|=\sup\{\|xy\|:y\in V, \|y\|\leq 1\}$.
\end{lem}

\noindent
\textbf{($\bigstar$)} Let $0<p\leq 1$, $X$ be a set and $\mathcal{A}$ be a unital Banach algebra with unit $\mathbf e$. Let $(F,|\cdot|)$ denote a commutative $p$-Banach algebra of complex valued functions on $X$ with pointwise multiplication, containing $\mathbf1(x)=1$ for all $x\in X$ as its unit. Let $\mathcal{F}$ denote a family of functions from $X$ to $\mathcal{A}$  with the following properties:
\begin{enumerate}
\item $\mathcal{F}$ is an algebra with pointwise multiplication.
\item If $a_1,a_2,\ldots,a_n \in \mathcal{A}$ and $f_1,f_2,\ldots,f_n \in F$, then $a_1f_1+a_2f_2+\cdots+a_nf_n \in \mathcal{F}.$
\item $\mathcal{F}$ is a $p$-Banach algebra with $p$-norm $\|\cdot\|_{\mathcal F}$ and $\|af\|_{\mathcal F} = \|a\|^p |f|$ for all $a\in \mathcal A$ and $f \in F$.
\item The set $\{a_1f_1+a_2f_2+\cdots+a_nf_n:a_i \in \mathcal A, f_i \in F, n\in \mathbb N\}$ is dense in $\mathcal{F}$.
\item If $x=a_1f_1+a_2f_2+\cdots+a_nf_n$, then for every multiplicative linear functional $M$ of $F$,
		$$\|a_1M(f_1)+a_2M(f_2)+\cdots+a_nM(f_n)\| \leq \|x\|_{\mathcal{F}}^{1/p}.$$
\end{enumerate}
It follows from (iv) and (v) that every $M$ gives rise to a continuous homomorphism $\bold{M}$ from $\mathcal{F}$ to $\mathcal{A}$ (known as \emph{generated homomorphism}), with the property $\bold{M}(af)=aM(f)$ for all $a\in\mathcal{A}$ and $f\in F$.

Following is a version of \cite[Theorem 3]{Bo} by Bochner and Phillips  for $p-$Banach algebras. The case for $p=1$ was proved by them.

\begin{thm}\label{th1}
Let $F, \mathcal F$ and $\mathcal A$ be as described in \textbf{($\bigstar$)}, and let $x \in \mathcal{F}$. If $\bold{M}(x)$ has a left inverse in $\mathcal{A}$ for every generated homomorphism $\bold{M}$, then $x$ has a left inverse in $\mathcal{F}$. If $F$ has no unit, then the element $\lambda \mathbf e\cdot\mathbf1+x\in\mathcal{F}_\mathbf1$, the unitization of $\mathcal{F}$, has a left inverse of the form $\lambda' \mathbf e\cdot\mathbf1+y$ whenever $\lambda \mathbf e+\bold{M}(x)$ has a left inverse in $\mathcal{A}$ for every generated homomorphism $\bold{M}$.
\end{thm}
\begin{proof}
We take $\mathcal{B}$ to be $\mathcal{F}$ in Lemma \ref{lbo10} and consider $I$ to be an arbitrary maximal ideal of $\mathcal{F}$ and $L_f$ to be an element of $B_0$ for $f\in\mathcal{F}$. Notice that $\mathbf ef\in\mathcal{F}$ commutes with each $ag\in\mathcal{F}$ for all $a\in\mathcal{A}$ and $f,g\in F$ and hence by (iv) with all elements of $\mathcal{F}$. So, by Lemmas \ref{lbo9}, \ref{lbo10} and \ref{lbo13}, $L_{\mathbf ef}$ can be seen as a ring homomorphism from $F$ to $\mathbb{C}$. This means that there is $M(f)$ such that $$L_{\mathbf ef}=\mathbf{e_1}M(f) \quad (\mathbf{e_1}=L_\mathbf{e1} \ \text{is the identity operator}).$$ Let $x=a_1f_1+a_2f_2+\dots +a_nf_n\in\mathcal{F}$. Then \begin{align*} L_x &=L_{(a_1f_1+a_2f_2+\dots +a_nf_n)} \\ &=L_{a_1 \mathbf1}L_{\mathbf{e}f_1}+L_{a_2 \mathbf1}L_{\mathbf{e}f_2}+ \dots +L_{a_n \mathbf1}L_{\mathbf{e}f_n} \\ &=L_{a_1 \mathbf1}M(f_1)+L_{a_2 \mathbf1}M(f_2)+\dots +L_{a_n \mathbf1}M(f_n) \\ &=L_{(a_1M(f_1)+a_2M(f_2)+\dots +a_nM(f_n))\mathbf1} \\ &=L_{\mathbf{M}(x) \mathbf1}. \end{align*} By Lemma \ref{lbo13} and the property (v), this relation is valid for all $x\in\mathcal{F}$. If $y$ is a left inverse of $M(x)$ in $\mathcal{A}$, then $L_yL_x=L_{yx}=L_{\mathbf M(yx) \mathbf1}=L_{yM(x) \mathbf1}=L_\mathbf{e1}=\mathbf e_1$. Thus $L_x$ has a left inverse, and by Lemma \ref{lbo10}, $x$ has a left inverse in $\mathcal{F}$.

If $F$ has no unit, then we take unitisation of $F$ and get the desired result.
\end{proof}

Let $0<p\leq1$, $\omega$ be a weight on $\mathbb Z$, and let $\mathcal A$ be a unital Banach algebra. Then the $p$-Banach space $$\ell^p(\mathbb Z,\omega,\mathcal A)=\{f:\mathbb Z \to \mathcal A:|f|_{p\omega}=\sum_{n \in \mathbb Z}\|f(n)\|^p\omega(n)^p<\infty\}$$ is a unital $p$-Banach algebra with the $p$-norm $|\cdot|_{p\omega}$ and the convolution product
\begin{eqnarray}\label{CONV}
(f\star g)(n)=\sum_{m\in \mathbb Z}f(m)g(n-m)\quad(n\in \mathbb Z) \ (f,g \in \ell^p(\mathbb Z,\omega,\mathcal A)).
\end{eqnarray}

\begin{lem}\label{lmcz}
Let $0<p\leq 1$, $\omega$ be a weight on $\mathbb{Z}$, $\mathcal{A}$ be a unital Banach algebra, $F=\{f\in C(\mathbb{T},\mathbb{C}): \widehat{f}\in\ell^p(\mathbb{Z},\omega,\mathbb{C})\}$, and let $\mathcal{F}=\{f\in C(\mathbb{T},\mathcal{A}): \widehat{f}\in\ell^p(\mathbb{Z},\omega,\mathcal{A})\}$. Then $F$ and $\mathcal{F}$ satisfy the conditions (i) to (v) in \textbf{($\bigstar$)}.
\end{lem}
\begin{proof}
The conditions (i) to (iv) are clearly satisfied. To check condition (v), let $\phi\in\Delta(F)$, and let $x=a_1f_1+a_2f_2+\dots+a_nf_n\in\mathcal{F}$ for some $n\in\mathbb{N}$ where $a_i\in\mathcal{A}$ and $f_i\in F$ for $1\leq i\leq n$. Since the Gel'fand space $\Delta(F)$ of $F$ is identified with $\Gamma(\rho_{1,\omega},\rho_{2,\omega})$ via the map $\phi\mapsto\phi_z$, where $\phi_z(f)=f(z) \ (f\in F)$, there is some $z\in\Gamma(\rho_{1,\omega},\rho_{2,\omega})$ such that $\phi=\phi_z$. Since $x\in\mathcal{F}$, $$\|x\|_{\mathcal{F}}=\sum_{k\in\mathbb{Z}}\|a_1\widehat{f_1}(k)+a_2\widehat{f_2}(k)+\dots+a_n\widehat{f_n}(k)\|^p\omega(k)^p<\infty.$$ So, $c=\left(\|a_1\widehat{f_1}(k)+a_2\widehat{f_2}(k)+\dots+a_n\widehat{f_n}(k)\|^p\omega(k)^p\right)_{k\in\mathbb{Z}}\in \ell^1(\mathbb{Z})\subset \ell^{1/p}(\mathbb{Z})$ and thus \begin{align*} \|c\|_{1/p}&=\left(\sum_{k\in\mathbb{Z}}(\|a_1\widehat{f_1}(k)+a_2\widehat{f_2}(k)+\dots+a_n\widehat{f_n}(k)\|)\omega(k)\right)^p \\ &\leq \sum_{k\in\mathbb{Z}}\|a_1\widehat{f_1}(k)+a_2\widehat{f_2}(k)+\dots+a_n\widehat{f_n}(k)\|^p\omega(k)^p =\|c\|_1. \end{align*} Let $z\in\Gamma(\rho_{1,\omega},\rho_{2,\omega})$, i.e., $\sup_{k\in\mathbb{N}}\omega(-k)^{-1/k} \leq|z|\leq \inf_{k\in\mathbb{N}}\omega(k)^{1/k}$. Then $\omega(-k)^{-1/k}\leq|z|\leq\omega(k)^{1/k}$ for all $k\in\mathbb{N}$ and $|z|^0=1\leq\omega(0)$. This implies that $|z|^k\leq\omega(k)$ for all $k\in\mathbb{Z}$. Thus,
\begin{eqnarray*}
&& \|a_1\phi_z(f_1)+a_2\phi_z(f_2)+\dots+a_n\phi_z(f_n)\|\\
&=& \|a_1f_1(z)+a_2f_2(z)+\dots+a_nf_n(z)\|\\
&=&\|\sum_{k\in\mathbb{Z}} (a_1\widehat{f_1}(k)+a_2\widehat{f_2}(k)+\dots+a_n\widehat{f_n}(k))z^k\|\\
&\leq & \sum_{k\in\mathbb{Z}} \|(a_1\widehat{f_1}(k)+a_2\widehat{f_2}(k)+\dots+a_n\widehat{f_n}(k))\|\omega(k)\\ &\leq & \left(\sum_{k\in\mathbb{Z}} \|(a_1\widehat{f_1}(k)+a_2\widehat{f_2}(k)+\dots+a_n\widehat{f_n}(k))\|^p \omega(k)^p\right)^{1/p}\\
&=& \|x\|_\mathcal{F}^{1/p}.
\end{eqnarray*}
\end{proof}

\begin{proof} [\textbf{Proof of Theorem \ref{thm1}}]
Let $\mathcal{A}_p(\omega)=\{g \in C(\mathbb{T},\mathcal{A}):\widehat{g} \in \ell^p(\mathbb{Z},\omega,\mathcal{A})\}$. Then $\mathcal{A}_p(\omega)$ is a unital $p$-Banach algebra with pointwise operations and $p$-norm $\|g\|_{p\omega}=|\widehat{g}|_{p\omega} \; (g \in \mathcal{A}_p(\omega))$. Also, $g \in C(\mathbb{T},\mathcal{A})$ has $p\omega$-ACFS if and only if $g \in \mathcal{A}_p(\omega)$ if and only if $\widehat{g}\in \ell^p(\mathbb{Z},\omega,\mathcal{A})$. Since $f$ has $p\omega$-ACFS, $f(z)=\sum_{n\in \mathbb Z}\widehat f(n) z^n\;(z \in \mathbb T)$ and $\|f\|_{p\omega}=\sum_{n\in \mathbb Z}\|\widehat f(n)\|^p\omega(n)^p<\infty$. Let $z\in\Gamma(\rho_{1,\omega},\rho_{2,\omega})$. Then $|z|^n\leq\omega(n)$ for all $n\in\mathbb{Z}$. Since $\sum_{n\in \mathbb Z}\|\widehat f(n)\|^p\omega(n)^p<\infty$, there is $n_0\in\mathbb{N}$ such that $\sum_{|n|\geq n_0}\|\widehat f(n)\|^p\omega(n)^p<1$. This implies that $$\sum_{n\in\mathbb{Z}}\|\widehat f(n)\||z|^n\leq\sum_{|n|<n_0}\|\widehat f(n)\||z|^n + \sum_{|n|\geq n_0}\|\widehat f(n)\|^p\omega(n)^p<\infty.$$ So, $\sum_{n\in\mathbb{Z}}\widehat f(n)z^n$ converges for all $z\in\Gamma(\rho_{1,\omega},\rho_{2,\omega})$.

Define $\Psi:\Gamma(\rho_{1,\omega},\rho_{2,\omega}) \to \mathcal{A}$ by $$\Psi(z)=\sum_{n \in \mathbb{Z}} \widehat f(n)z^n \quad (z \in \Gamma (\rho_{1,\omega},\rho_{2,\omega})).$$ Then $\Psi$ is continuous and $\Psi(z)=f(z)$ for all $z \in \mathbb T$. Take any $z_0 \in \mathbb{T}$. Since $\Psi(z_0)$ is left invertible in $\mathcal{A}$ and the set $U_l$ of all left invertible elements in $\mathcal{A}$ is open, there is $\epsilon>0$ such that the open ball $B(\Psi(z_0),\epsilon)$ is contained in $U_l$. By continuity of $\Psi$, there is some $\delta_{z_0}>0$ such that $\Psi(t)\in B(\Psi(z_0),\epsilon)$ for all $t \in B(z_0,\delta_{z_0})\cap \Gamma(\rho_{1,\omega},\rho_{2,\omega})$, i.e., $\Psi(t)\in U_l$ for all $t \in B(z_0,\delta_{z_0})\cap \Gamma(\rho_{1,\omega},\rho_{2,\omega})$. By compactness of $\mathbb{T}$, we get $r_1,r_2\in\mathbb{R}$ such that $\rho_{1,\omega} \leq r_1 \leq 1 \leq r_2 \leq \rho_{2,\omega}$ and $\Psi(t)\in U_l$ for all $t \in \Gamma(r_1,r_2)$. So, $f(z)=\sum_{n \in \mathbb{Z}} \widehat f(n) z^n \in U_l$ for all $z \in \Gamma(r_1,r_2)$. Define $\nu : \mathbb{Z} \to[1,\infty)$ as follows: If $\rho_{1,\omega}=\rho_{2,\omega}$, then take $\nu = \omega$; otherwise define $$\nu(n)=\begin{cases} r_1^n & n\leq 0 \\ r_2^n & n\geq 0\end{cases}.$$ Then $\nu$ is a weight on $\mathbb{Z}$ which is constant if and only if $\omega$ is constant and $1 \leq \nu \leq \omega$. This gives $\mathcal{A}_p(\omega) \subset \mathcal{A}_p(\nu)$ and so $f \in \mathcal{A}_p(\nu)$. Let $F=\{g\in C(\mathbb{T},\mathbb{C}):\widehat{g}\in \ell^p(\mathbb{Z},\nu) \}$, and let $\mathcal{F}= \mathcal{A}_p(\nu)$. Then $f \in \mathcal{F}$. The Gel'fand space $\Delta(F)$ of $F$ is identified with $\Gamma(r_1,r_2)$ via the map $z \mapsto \phi_z$, where $\phi_z(g)=g(z) \ (g \in F)$. For every generated homomorphism $\Phi_z$, $\Phi_z(f)=f(z)$ is left invertible in $\mathcal{A}$. So, by Lemma \ref{lmcz} with Theorem \ref{th1}, $f$ has left inverse in $\mathcal{F}$, say, $g$. Then $g$ is our required function in $\mathcal{A}_p(\nu)$. This completes the proof.
\end{proof}

\section{Proof of Theorem  \ref{thm2}}
Let $1<p<\infty$, and let $q$ be its conjugate index. Let $\omega$ be a weight on $\mathbb Z$ satisfying $\omega^{-q}\star\omega^{-q}\leq\omega^{-q}$, and let $\mathcal A$ be a unital Banach algebra. Then, by  \cite{SG}, $$\ell^p(\mathbb Z,\omega,\mathcal A)=\{f:\mathbb Z \to \mathcal A:|f|_{p\omega}=\left(\sum_{n \in \mathbb Z}\|f(n)\|^p\omega(n)^p\right)^\frac{1}{p}<\infty\}$$ is a unital Banach algebra with the norm $|\cdot|_{p\omega}$ and the convolution product defined in the equation (\ref{CONV}).

\begin{lem}\label{lm*}
Let $1<q<\infty,\ 0<\rho_1<1<\rho_2$, and let $\omega$ be a weight on $\mathbb{Z}$ satisfying $\sum_{n\in\mathbb{Z}}\omega(n)^{-q}<\infty$ and $\omega^{-q}\star \omega^{-q}\leq\omega^{-q}$. If $\sum_{n\in\mathbb{Z}}x^{qn}\omega(n)^{-q}<\infty$ for all $x\in(\rho_1,\rho_2)$, then $\sum_{n\in\mathbb{Z}}x^{qn}\omega(n)^{-q}\leq 1$ for $x\in[\rho_1,\rho_2]$.
\end{lem}
\begin{proof}
Define $f:(\rho_1^q,\rho_2^q)\to(0,\infty)$ by $f(x)=\sum_{n\in\mathbb{Z}}x^n\omega(n)^{-q} \quad(x \in (\rho_1^q,\rho_2^q)).$ Let $x\in(\rho_1^q,\rho_2^q)$. The fact $\omega^{-q}\star\omega^{-q}\leq\omega^{-q}$ implies that
\begin{align}\label{EQ}
f(x) =\sum_{n\in\mathbb{Z}}x^n\omega(n)^{-q}&\geq \sum_{n\in\mathbb{Z}}x^n(\omega^{-q}\star \omega^{-q})(n) \nonumber \\
&=\sum_{n\in\mathbb{Z}} x^n \left(\sum_{k\in\mathbb{Z}}\omega(k)^{-q}\omega(n-k)^{-q}\right) \nonumber \\ &=\sum_{k\in\mathbb{Z}} x^{k}\omega(k)^{-q} \left(\sum_{n\in\mathbb{Z}} x^{n-k}\omega(n-k)^{-q} \right) \nonumber\\&
=f(x)^2.
\end{align}
Thus $0\leq f(x)\leq 1$. Define $g:[1,\rho_2^q)\to(0,\infty)$ by $g(x)=\sum_{n\in\mathbb{N}_0} x^n\omega(n)^{-q} \quad(x\in [1,\rho_2^q))$. Since the radius of convergence of $\sum_{n\in\mathbb{N}_0} x^n\omega(n)^{-q}$ is at least $\rho_2^q$, the function $g$ is continuous on $[1,\rho_2^q)$. Also, if $x\in[1,\rho_2^q)$, then $0\leq g(x)\leq f(x) \leq1$. Moreover, if $1\leq x<y<\rho_2^q$, then $g(x)<g(y)$. So, $g$ is continuous, bounded and increasing on $[1,\rho_2^q)$ and therefore $g$ is uniformly continuous and can be extended as a continuous function on $[1,\rho_2^q]$. This gives $g(\rho_2^q)=\sum_{n\in\mathbb{N}_0} \rho_2^{qn}\omega(n)^{-q}<\infty$. Now,
\begin{align*} \sum_{n\in\mathbb{Z}} \rho_2^{nq}\omega(n)^{-q} &=\sum_{n\in\mathbb{N}} \rho_2^{-nq}\omega(-n)^{-q} + \sum_{n\in\mathbb{N}_0} \rho_2^{nq}\omega(n)^{-q} \\ &\leq \sum_{n\in\mathbb{N}} \omega(-n)^{-q} + \sum_{n\in\mathbb{N}_0} \rho_2^{nq}\omega(n)^{-q} \\ &\leq \sum_{n\in\mathbb{Z}} \omega(n)^{-q} + \sum_{n\in\mathbb{N}_0} \rho_2^{nq}\omega(n)^{-q} \\ &<\infty. \end{align*} Similarly it can be shown that $\sum_{n\in\mathbb{Z}} \rho_1^{nq}\omega(n)^{-q}<\infty$. Using the same computation as in inequality (\ref{EQ}) we have $\sum_{n\in\mathbb{Z}}\rho_1^{qn}\omega(n)^{-q}\leq 1$ and $\sum_{n\in\mathbb{Z}}\rho_2^{qn}\omega(n)^{-q}\leq 1$.
\end{proof}

\begin{prop}\label{prz}
Let $1<p<\infty$, $q$ be its conjugate index, and let $\omega$ be a weight on $\mathbb{Z}$ satisfying $\sum_{n\in\mathbb{Z}}\omega(n)^{-q}<\infty$ and $\omega^{-q}\star \omega^{-q}\leq\omega^{-q}$. Then the Gel'fand space $\Delta(\ell^p(\mathbb{Z},\omega,\mathbb{C}))$ of $\ell^p(\mathbb{Z},\omega,\mathbb{C})$ is the closed annulus $\Gamma_\omega=\Gamma(\rho_{1,\omega},\rho_{2,\omega})$.
\end{prop}
\begin{proof}
Let $f=(a_n)_{n\in\mathbb{Z}}\in F=\ell^p(\mathbb{Z},\omega,\mathbb{C})$. Then $\|f\|_{p\omega}=\left(\sum_{n \in \mathbb Z}|a_n|^p\omega(n)^p\right)^{1/p}<\infty$. If $z\in\mathbb{T}$, then \begin{align*}\sum_{n\in\mathbb{Z}}|a_n||z|^n=\sum_{n\in\mathbb{Z}}|a_n|\omega(n)\omega(n)^{-1}\leq \|f\|_{p\omega}\left(\sum_{n\in\mathbb{Z}}\omega(n)^{-q}\right)^{1/q}<\infty. \end{align*} If $1<|z|<\rho_{2,\omega}$, then $\frac{|z|}{\rho_{2,\omega}}<1$ and so, \begin{align*}\sum_{n\in\mathbb{Z}}|a_n||z|^n&=\sum_{n\in\mathbb{N}}|a_{-n}||z|^{-n} + \sum_{n\in\mathbb{N}_0}|a_n|\omega(n)\frac{|z|^n}{\omega(n)}\\ &\leq\sum_{n\in\mathbb{N}}|a_{-n}|\omega(-n)\omega(-n)^{-1} + \sum_{n\in\mathbb{N}_0}|a_n|\omega(n)\frac{|z|^n}{\rho_{2,\omega}^n}\\ &\leq \|f\|_{p\omega}\left(\sum_{n\in\mathbb{Z}}\omega(n)^{-q}\right)^{1/q} + \left(\sum_{n\in\mathbb{N}_0}|a_n|^p\omega(n)^p\right)^{1/p} \left(\sum_{n\in\mathbb{N}_0}\frac{|z|^{nq}}{\rho_{2,\omega}^{nq}}\right)^{1/q}\\&<\infty. \end{align*} If $\rho_{1,\omega}<|z|<1$, then $\sum_{n\in\mathbb{Z}}|a_n||z|^n<\infty$ can be shown similarly. By Lemma \ref{lm*}, the sequences $(\rho_{1,\omega}^n)_{n\in\mathbb{Z}}$ and $(\rho_{2,\omega}^n)_{n\in\mathbb{Z}}$ are in $\ell^q(\mathbb{Z},\omega^{-1})$ and so if $|z|=\rho_{1,\omega}$ or $|z|=\rho_{2,\omega}$, then
\begin{align*}
\sum_{n\in\mathbb{Z}}|a_n||z|^n  \leq \sum_{n\in\mathbb{Z}}|a_n||z|^n \frac{\omega(n)}{\omega(n)} \leq \left(\sum_{n\in\mathbb{Z}} |a_n|^p \omega(n)^p \right)^{1/p} \left(\sum_{n\in\mathbb{Z}} |z|^{qn}  \omega(n)^{-q}\right)^{1/q}<\infty.
\end{align*}
Thus, $\sum_{n\in\mathbb{Z}}|a_n| |z|^n<\infty$ for all $z\in\Gamma_\omega$.

Define $\phi_z:F\to\mathbb{C}$ by $$\phi_z(f)=\sum_{n\in\mathbb Z}a_nz^n \quad (f=(a_n)\in F).$$ Then it is clear that $\phi_z\in\Delta(F)$. Let $f=(a_n)\in F$. Then we may write $f$ as $f=\sum_{n\in\mathbb{Z}}a_n\delta_n=\sum_{n\in\mathbb{Z}}a_n\delta_1^n$, where $\delta_n(n)=1$ and $0$ otherwise. Let $\phi\in\Delta(F)$. Let $z=\phi(\delta_1)$. Then for $f\in F$, $\phi(f)=\sum_{n\in\mathbb{Z}}a_n\phi(\delta_1)^n=\sum_{n\in\mathbb{Z}}a_nz^n$. Since $\delta_1$ is invertible, $z\neq0$. Since $\phi$ is continuous, $|z|^n=|\phi(\delta_n)|\leq\|\delta_n\|\leq\omega(n)$ for all $n\in\mathbb{Z}$. Therefore $z\in\Gamma_\omega$ and $\phi=\phi_z$. Hence, $\Delta(F)=\{\phi_z:z\in\Gamma_\omega\}$.
\end{proof}

\begin{lem}\label{lmcz2}
Let $1<p<\infty, q$ be the conjugate index of $p$, $\omega$ be a weight on $\mathbb{Z}$ satisfying $\sum_{n\in\mathbb{Z}}\omega(n)^{-q}<\infty$ and $\omega^{-q}\star\omega^{-q}\leq\omega^{-q}$, $\mathcal{A}$ be a unital Banach algebra, $F=\{f\in C(\mathbb{T},\mathbb{C}): \widehat{f}\in\ell^p(\mathbb{Z},\omega,\mathbb{C})\}$, and let $\mathcal{F}=\{f\in C(\mathbb{T},\mathcal{A}): \widehat{f}\in\ell^p(\mathbb{Z},\omega,\mathcal{A})\}$. Then $F$ and $\mathcal{F}$ satisfy the conditions (i) to (v) in \textbf{($\bigstar$)}.
\end{lem}
\begin{proof}
It is clear that the conditions (i) to (iv) are satisfied. For condition (v), let $\phi\in\Delta(F)$, and let $x=a_1f_1+a_2f_2+\dots+a_nf_n\in\mathcal{F}$ for some $n\in\mathbb{N}$ where $a_i\in\mathcal{A}$ and $f_i\in F$ for $1\leq i\leq n$. By Proposition \ref{prz}, $\phi=\phi_z$ for some $z\in\Gamma(\rho_{1,\omega},\rho_{2,\omega})$. Then
\begin{align*} &\|a_1\phi_z(f_1)+a_2\phi_z(f_2)+\dots+a_n\phi_z(f_n)\| \\&= \|a_1f_1(z)+a_2f_2(z)+\dots+a_nf_n(z)\| \\ &= \|\sum_{k\in\mathbb{Z}} (a_1\widehat{f_1}(k)+a_2\widehat{f_2}(k)+\dots+a_n\widehat{f_n}(k))z^k\| \\ &\leq \sum_{k\in\mathbb{Z}} \|(a_1\widehat{f_1}(k)+a_2\widehat{f_2}(k)+\dots+a_n\widehat{f_n}(k))\|\omega(k) |z|^k \omega(k)^{-1} \\ &\leq \left(\sum_{k\in\mathbb{Z}} \|(a_1\widehat{f_1}(k)+a_2\widehat{f_2}(k)+\dots+a_n\widehat{f_n}(k))\|^p \omega(k)^p\right)^{1/p} \left(\sum_{k\in\mathbb{Z}}|z|^{kq}\omega(k)^{-q}\right)^{1/q} \\&\leq \|x\|_\mathcal{F}, \end{align*} as $\sum_{k\in\mathbb{Z}}|z|^{kq}\omega(k)^{-q} \leq1$ by Lemma \ref{lm*}.
\end{proof}

We shall require the following two lemmas to prove Theorem \ref{thm2}.

\begin{lem} \cite[Lemma 1]{Da} \label{lm1}
Let $\omega$ be a weight on $\mathbb{Z}$.
\begin{enumerate}
\item If $\rho_{2,\omega}>1$ and $\gamma \in (0,1)$, then there is a positive constant $K$ such that $K\omega(n)\geq e^{n^\gamma}$ for all $n \in \mathbb{N}$.
\item If $\rho_{1,\omega}<1$ and $\gamma \in (0,1)$, then there is a positive constant $K$ such that $K\omega(-n)\geq e^{n^\gamma}$ for all $n \in \mathbb{N}$.
\item If $\rho_{1,\omega}<1<\rho_{2,\omega}$ and $\gamma \in (0,1)$, then there is a positive constant $K$ such that $K\omega(n)\geq e^{|n|^\gamma}$ for all $n \in \mathbb{Z}$.
\end{enumerate}
\end{lem}

\begin{lem} \cite[Lemma 2]{Da} \label{lm2}
Let $\omega$ be an almost monotone algebra weight on $\mathbb{Z}$ such that $\rho_{2,\omega}=1$. Define $\widetilde{\omega}:\mathbb{Z} \to [1,\infty)$ by $\widetilde{\omega}(n)=\max \{\omega(k):0\leq k \leq n \}$ if $n \in \mathbb{N}_0$ and $\widetilde{\omega}(n)=\omega(n)$ if $n<0$. Then $\widetilde{\omega}$ is an almost monotone algebra weight on $\mathbb{Z}$ and there is $K>0$ such that $\widetilde{\omega}(n)\leq K \omega(n)$ for all $n \in \mathbb{Z}$.
\end{lem}

\begin{proof} [\textbf{Proof of Theorem \ref{thm2}}]
Let $\mathcal{A}_p(\omega)=\{g \in C(\mathbb{T},\mathcal{A}):\widehat{g} \in \ell^p(\mathbb{Z},\omega,\mathcal{A})\}$. Then $\mathcal{A}_p(\omega)$ is a unital Banach algebra with pointwise operations and norm $\|g\|_{p\omega}=|\widehat{g}|_{p\omega} \; (g \in \mathcal{A}_p(\omega))$. Also, $g \in C(\mathbb{T},\mathcal{A})$ has $p\omega$-ACFS if and only if $g \in \mathcal{A}_p(\omega)$ if and only if $\widehat{g}\in \ell^p(\mathbb{Z},\omega,\mathcal{A})$. Since $f$ has $p\omega$-ACFS, $f(z)=\sum_{n\in \mathbb Z}a_n z^n\;(z \in \mathbb T)$, where $a_n \in \mathcal A$ are the Fourier coefficients of $f$, and $\|f\|_{p\omega}=\left(\sum_{n \in \mathbb Z}\|a_n\|^p\omega(n)^p\right)^{1/p}<\infty$. As done in Proposition \ref{prz}, taking $\|a_n\|$ in place of $|a_n|$, we get that $\sum_{n\in\mathbb{Z}}a_n z^n$ converges for all $z\in\Gamma(\rho_{1,\omega},\rho_{2,\omega})$.

Define $\Psi:\Gamma(\rho_{1,\omega},\rho_{2,\omega})\to \mathcal{A}$ by $$\Psi(z)=\sum_{n \in \mathbb{Z}} a_nz^n \quad (z \in \Gamma(\rho_{1,\omega},\rho_{2,\omega})).$$ Then $\Psi$ is continuous and $\Psi(z)=f(z)$ for all $z \in \mathbb T$. As done in Theorem \ref{thm1}, we get $r_1,r_2\in\mathbb{R}$ such that $\rho_{1,\omega} \leq r_1 \leq 1 \leq r_2 \leq \rho_{2,\omega}$ and $\Psi(t)$ is left invertible for all $t \in \Gamma(r_1,r_2)$. So, $f(z)=\sum_{n \in \mathbb{Z}} a_n z^n$ is left invertible in $\mathcal{A}$ for all $z \in \Gamma(r_1,r_2)$. Observe that if $\rho_{1,\omega}<1$, then $\rho_{1,\omega}<r_1<1$; if $\rho_{1,\omega}=1$, then $r_1=1$; if $\rho_{2,\omega}>1$, then $1<r_2<\rho_{2,\omega}$; and if $\rho_{2,\omega}=1$, then $r_2=1$. Define $\nu : \mathbb{Z} \to [1,\infty)$ as follows:\\
If $\rho_{1,\omega} = \rho_{2,\omega}$, then take $\nu = \omega$. \\
If $\rho_{1,\omega}=1$ and $\rho_{2,\omega}>1$, then, by lemma \ref{lm2}, there is $K_1>0$ such that $\omega(n)\leq \widetilde{\omega}(n) \leq K_1\omega(n)$ for all $n<0$, where $\widetilde{\omega}(n)=\max\{\omega(k):n \leq k \leq -1\}$ for all $n<0$. For $\gamma \in (0,1)$, we get $K_2>0$ such that $e^{|n|^\gamma} \leq K_2\omega(n)$ for all $n \in \mathbb{N}_0$. Set $K=\max\{K_1,K_2\}$ and define $$\nu(n)=\begin{cases} \widetilde{\omega}(n) & n<0 \\ r_2^\frac{n}{2} e^\frac{|n|^\gamma}{2} & n\geq 0\end{cases}.$$
Let $\rho_{1,\omega}<1$ and $\rho_{2,\omega}=1$. If $\widetilde{\omega}(n)=\max\{\omega(k):0\leq k \leq n\}$ for all $n\in \mathbb{N}_0$, then, by lemma \ref{lm2}, there is $K_1>0$ such that $\omega(n)\leq \widetilde{\omega}(n) \leq K_1\omega(n)$ for all $n\in \mathbb{N}_0$. Take $\gamma \in (0,1)$. Then there is $K_2>0$ such that $e^{|n|^\gamma} \leq K_2\omega(n)$ for all $n<0$. Let $K=\max\{K_1,K_2\}$ and define $$\nu(n)=\begin{cases} r_1^\frac{n}{2} e^\frac{|n|^\gamma}{2} & n<0 \\ \widetilde{\omega}(n) & n\geq 0\end{cases}.$$
If $\rho_{1,\omega}<1<\rho_{2,\omega}$, then for $\gamma \in (0,1)$, we get $K>0$ such that $e^{|n|^\gamma} \leq K\omega(n)$ for all $n\in\mathbb{Z}$. Define $$\nu(n)=\begin{cases} r_1^\frac{n}{2} e^\frac{|n|^\gamma}{2} & n<0 \\ r_2^\frac{n}{2} e^\frac{|n|^\gamma}{2} & n\geq 0\end{cases}.$$
Then $\nu$ is an almost monotone algebra weight on $\mathbb{Z}$ such that $1 \leq \nu \leq K\omega$ for some $K>0$. This gives $\mathcal{A}_p(\omega) \subset \mathcal{A}_p(\nu)$ and so $f \in \mathcal{A}_p(\nu)$. Let $F=\{g\in C(\mathbb{T},\mathbb{C}):\widehat{g}\in \ell^p(\mathbb{Z},\nu) \}$, and let $\mathcal{F}= \mathcal{A}_p(\nu)$. Then $f \in \mathcal{F}$. The Gel'fand space $\Delta(F)$ of $F$ is $\Gamma(\sqrt{r_1},\sqrt{r_2})$ via $z \mapsto \phi_z$, where $\phi_z(g)=g(z)\ (g \in F)$. Note that any generated homomorphism on $\mathcal{F}$ has the form $\Phi_z(g)=g(z) \ (g\in\mathcal{F})$ for some $z\in\Gamma(\sqrt{r_1},\sqrt{r_2})$. As $\Phi_z(f)=f(z)$ is left invertible in $\mathcal{A}$ for all $z\in\Gamma(\sqrt{r_1},\sqrt{r_2})$, by Lemma \ref{lmcz} and Theorem \ref{th1}, $f$ has left inverse in $\mathcal{F}$, say, $g$. Then $g$ is our required function in $\mathcal{A}_p(\nu)$. This completes the proof.
\end{proof}

We write important corollaries of Theorems \ref{thm1} and \ref{thm2} which would be helpful in the coming results.

\begin{cor} \label{crm1}
Let $0<p\leq1$. Let $\omega$ be a weight on $\mathbb{Z}$, $\mathcal A$ be a unital Banach algebra, and let $f \in \ell^p(\mathbb{Z},\omega,\mathcal A)$. If $\widehat{f}(z)=\sum_{n \in \mathbb{Z}} f(n) z^n$ is invertible in $\mathcal A$ for all $z \in \mathbb{T}$, then there is a weight $\nu$ on $\mathbb{Z}$ such that $1 \leq \nu \leq \omega,\ \nu$ is constant if and only if $\omega$ is constant and $f^{-1} \in \ell^p(\mathbb{Z},\nu,\mathcal A)$. In particular, if $\omega$ is an admissible weight, then $f^{-1} \in \ell^p(\mathbb{Z},\omega,\mathcal A)$.
\end{cor}

\begin{cor} \label{crm2}
Let $1<p<\infty$. Let $\omega$ be an almost monotone algebra weight on $\mathbb{Z}$, $\mathcal A$ be a unital Banach algebra, and let $f \in \ell^p(\mathbb{Z},\omega,\mathcal A)$.  If $\widehat{f}(z)=\sum_{n \in \mathbb{Z}} f(n) z^n$ is invertible in $\mathcal A$ for all $z \in \mathbb{T}$, then there is an almost monotone algebra weight $\nu$ on $\mathbb{Z}$ such that $1 \leq \nu \leq K\omega$ for some $K>0$ and $f^{-1} \in \ell^p(\mathbb{Z},\nu,\mathcal A)$. In particular, if $\omega$ is an admissible weight, then $f^{-1} \in \ell^p(\mathbb{Z},\omega,\mathcal A)$.
\end{cor}

\section{Proof of Theorem \ref{thr}}

Let $\mathcal{A}$ be a unital Banach algebra, and let $\omega$ be a weight on $\mathbb{R}$, i.e., $\omega:\mathbb R \to[1,\infty)$ is Borel measurable and $\omega(x+y)\leq \omega(x)\omega(y)$ for all $x,y \in \mathbb R$. Then $$L^1(\mathbb{R},\omega,\mathcal{A})=\{f:\mathbb{R}\to\mathcal{A}:\|f\|_{\omega}=\int_\mathbb{R}\|f(x)\|\omega(x)dx<\infty\}$$ is a Banach algebra with the norm $\|\cdot\|_\omega$ and the convolution product given by $$(f\star g)(x)=\int_\mathbb{R}f(x-t)g(t)dt \quad (f,g\in L^1(\mathbb{R},\omega,\mathcal{A}), x\in \mathbb R).$$ Let $L^1(\mathbb{R},\omega,\mathcal{A})_\mathbf1$ denote the unitization of $L^1(\mathbb{R},\omega,\mathcal{A})$ with the unit $\mathbf1$.

Let $S_\omega=\{x+iy\in\mathbb{C}:\rho_{1,\omega}\leq x \leq\rho_{2,\omega}, y\in\mathbb{R}\}$. Let $a+ib\in S_\omega$, i.e., $\sup\{\log\omega(x)^{1/x}:x<0\} \leq a \leq \inf\{\log\omega(x)^{1/x}:x>0\}$. Then $e^{ax}\leq\omega(x)$ for all $x\in\mathbb{R}$. So, for $f\in L^1(\mathbb{R},\omega,\mathcal{A})$, $$\int_\mathbb{R} \|f(x)\||e^{(a+ib)x}|dx \leq \int_\mathbb{R} \|f(x)\|e^{ax}dx \leq \int_\mathbb{R} \|f(x)\|\omega(x)dx<\infty.$$ For $f\in L^1(\mathbb{R},\omega,\mathcal{A})$, the \emph{Fourier transform} $\widehat{f}$ of $f$ is given as $$\widehat{f}(z)=\int_\mathbb{R} f(x) e^{zx}dx \in \mathcal{A} \quad (z\in S_\omega).$$

\begin{lem}\label{lmcr}
Let $\omega$ be a weight on $\mathbb{R}, \mathcal{A}$ be a unital Banach algebra, $F=\{\widehat{f}:f\in L^1(\mathbb{R},\omega,\mathbb{C})\}$, and let $\mathcal{F}=\{\widehat{f}:f\in L^1(\mathbb{R},\omega,\mathcal{A})\}$. Then $F$ and $\mathcal{F}$ satisfy all conditions from (i) to (v) in \textbf{($\bigstar$)}.
\end{lem}
\begin{proof}
Clearly conditions (i) to (iv) are satisfied. For (v), let $\phi\in\Delta(F)$, and let $x=a_1f_1+a_2f_2+\dots+a_nf_n\in\mathcal{F}$ for some $n\in\mathbb{N}$ where $a_i\in\mathcal{A}$ and $f_i\in F$ for $1\leq i\leq n$. Since the Gel'fand space $\Delta(F)$ of $F$ is identified with $S_\omega$ via the map $\phi\mapsto\phi_z$, where $\phi_z(f)=\widehat f(z) \ (f\in F)$, there is some $z\in S_\omega$ such that $\phi=\phi_z$. So,
\begin{eqnarray*}
&&\|a_1\phi_z(f_1)+a_2\phi_z(f_2)+\dots+a_n\phi_z(f_n)\| \\
&= & \|a_1\widehat{f_1}(z)+a_2\widehat{f_2}(z)+\dots+a_n\widehat{f_n}(z)\|\\
&\leq & \int_{\mathbb{R}} \|a_1f_1(x)+a_2f_2(x)+\dots+a_nf_n(x)\||e^{zx}|dx\\
&\leq & \int_{\mathbb{R}} \|a_1f_1(x)+a_2f_2(x)+\dots+a_nf_n(x)\|\omega(x) dx \\&=&\|x\|_{\mathcal{F}}.
\end{eqnarray*}
\end{proof}

\begin{proof}[\textbf{Proof of Theorem \ref{thr}}]
Here $f\in L^1(\mathbb{R},\omega,\mathcal{A})$ and so by an application of the Riemann-Lebesgue lemma, $$\lim_{|z|\to\infty} \widehat{f}(z)=0.$$ Then there is $M>0$ such that $$\|\mathbf e-(\mathbf e+\widehat{f}(x+iy))\|=\|\widehat{f}(x+iy)\|<1$$ for all $x+iy\in S_\omega$ with $|y|>M.$ This means that $\mathbf e+\widehat{f}(x+iy)$ is left invertible in $\mathcal{A}$ for all $x+iy\in S_\omega$ with $|y|>M.$ Let $y\in\mathbb{R}$ with $|y|\leq M$. Then by hypothesis $\mathbf e+\widehat{f}(iy)$ is left invertible in $\mathcal{A}$. Since the set of all left invertible elements of $\mathcal{A}$ is open and $\mathbf e+\widehat{f}$ is continuous, there is $\delta_y>0$ such that $\mathbf e+\widehat{f}(z)$ is left invertible in $\mathcal{A}$ for all $z\in B(y,\delta_y) \cap S_\omega$. By compactness of $\{iy:|y|\leq M\}$, we get $r_1,r_2\in\mathbb{R}$ such that $\rho_{1,\omega}\leq r_1 \leq 0 \leq r_2 \leq \rho_{2,\omega}$ and $\mathbf e+\widehat{f}(z)$ is left invertible in $\mathcal{A}$ for all $z=x+iy\in [r_1,r_2] + i [-M,M]$. Define $\nu:\mathbb{R}\to [1,\infty)$ as follows:\\
If $\rho_{1,\omega}=0=\rho_{2,\omega}$, then take $\nu=\omega$.\\
If $\rho_{1,\omega}=0$ and $\rho_{2,\omega}>0$, then $$\nu(x)=\begin{cases} \omega(x) & x\leq0 \\ e^{r_2x} & x>0 \end{cases}.$$
If $\rho_{1,\omega}<0$ and $\rho_{2,\omega}=0$, then $$\nu(x)=\begin{cases} e^{r_1x} & x\leq0 \\ \omega(x) & x>0 \end{cases}.$$
If $\rho_{1,\omega}<0$ and $\rho_{2,\omega}>0$, then $$\nu(x)=\begin{cases} e^{r_1x} & x\leq0 \\ e^{r_2x} & x>0 \end{cases}.$$
Then $\nu$ is a weight on $\mathbb{R}$ which is constant if and only if $\omega$ is constant and $1\leq\nu\leq\omega$. This gives $L^1(\mathbb{R},\omega,\mathcal{A})\subset L^1(\mathbb{R},\nu,\mathcal{A})$ and so $f\in L^1(\mathbb{R},\nu,\mathcal{A})$. Let $F=\{\widehat{g}:g\in L^1(\mathbb{R},\nu,\mathbb{C})\}$, and let $\mathcal{F}=\{\widehat{g}:g\in L^1(\mathbb{R},\nu,\mathcal{A})\}$. Then $\widehat f\in \mathcal{F}$. The Gel'fand space $\Delta(F)$ of $F$ is the strip $S_\nu=\{x+iy:r_1\leq x \leq r_2,\ y\in \mathbb{R}\}$ identified by the mapping $z\mapsto\phi_z$, where $\phi_z(g)=g(z) \ (g\in F).$ For every generated homomorphism $\Phi_z$, $\mathbf e+\Phi_z(\widehat f)=\mathbf  e+\widehat f(z)$ is left invertible in $\mathcal{A}$. So, by Lemma \ref{lmcr} and Theorem \ref{th1}, $\mathbf  e+\widehat{f}$ has a left inverse in $\mathcal{F}$, say, $\mathbf  e+\widehat{g}$. This implies that $g\in L^1(\mathbb{R},\nu,\mathcal{A})$ is such that $\mathbf1+g$ is a left inverse of $\mathbf1+f$. This completes the proof.
\end{proof}

\section{Proof of Theorem \ref{thrp}}

Let $1<p<\infty$, and let $q$ be its conjugate index. Let $\mathcal{A}$ be a unital Banach algebra, and let $\omega$ be a weight on $\mathbb{R}$ satisfying $\omega^{-q}\star\omega^{-q}\leq\omega^{-q}$.
Then $$L^p(\mathbb{R},\omega,\mathcal{A})=\{f:\mathbb{R}\to\mathcal{A}:\|f\|_{p\omega}=\left(\int_\mathbb{R}\|f(x)\|^p \omega(x)^pdx\right)^\frac{1}{p}<\infty\}$$ is a Banach algebra with convolution as product and the norm $\|\cdot\|_{p\omega}$ \cite{SG}.

\begin{prop}
Let $1<p<\infty,\ q$ be its conjugate index, and let $\omega$ be a weight on $\mathbb{R}$ such that $\int_{\mathbb R}\omega^{-q}<\infty$ and $\omega^{-q}\star\omega^{-q}\leq\omega^{-q}$. Then the Gel'fand space of $L^p(\mathbb{R},\omega,\mathbb{C})$ is the strip $S_\omega$ and $\int_{\mathbb{R}}|e^{zx}|^q\omega(x)^{-q}dx\leq 1$ for $z \in S_\omega$.
\end{prop}
\begin{proof}
Let $f\in F=L^p(\mathbb{R},\omega,\mathbb{C})$, and let $z=a+ib\in S_\omega$. If $a=0$, then \begin{align*}\int_\mathbb{R}|f(x)||e^{zx}|dx =\int_\mathbb{R}|f(x)|dx &=\int_\mathbb{R}|f(x)|\omega(x)\omega(x)^{-1}dx\\ &\leq \left(\int_\mathbb{R}|f(x)|^p \omega(x)^p dx\right)^{1/p}\left(\int_\mathbb{R}\omega(x)^{-q}dx\right)^{1/q}\\ &<\infty.  \end{align*} If $\rho_{2,\omega}>0$ and $0<a<\rho_{2,\omega}$, then $a-\rho_{2,\omega}<0$ and so,
\begin{align*}\int_\mathbb{R}|f(x)||e^{zx}|dx &=\int_{-\infty}^0 |f(x)|e^{ax} dx + \int^\infty_0 |f(x)|e^{ax} \omega(x)\omega(x)^{-1}dx dx\\ & \leq \int_{-\infty}^0 |f(x)| dx + \int^\infty_0|f(x)|e^{ax} \omega(x) e^{-\rho_{2,\omega}x} dx\\ & \leq \int_{\mathbb R} |f(x)| dx + \left(\int^\infty_0 |f(x)|^p \omega(x)^p dx \right)^{1/p} \left(\int^\infty_0 e^{(a-\rho_{2,\omega})qx} dx\right)^{1/q}\\ &<\infty. \end{align*}
If $\rho_{1,\omega}<0$ and $\rho_{1,\omega}<a<0$, then $\int_\mathbb{R} |f(x)||e^{zx}|dx <\infty$ can be shown similarly.
If $\rho_{1,\omega}<0$ and/or $\rho_{2,\omega}>0$, then define $\zeta:(\rho_{1,\omega}q,\rho_{2,\omega}q)\to(0,\infty)$ by $$\zeta(y)=\int_{\mathbb{R}}e^{yx}\omega(x)^{-q}dx \quad (y\in(\rho_{1,\omega}q,\rho_{2,\omega}q)).$$ If $y\in(\rho_{1,\omega}q,\rho_{2,\omega}q)$, then
\begin{eqnarray}\label{EQQ}
\zeta(y)&=&\int_{\mathbb{R}}e^{yx}\omega(x)^{-q}dx \nonumber \\ &\geq& \int_{\mathbb{R}}e^{yx} \left(\int_{\mathbb{R}}\omega(x_1)^{-q} \omega(x-x_1)^{-q} dx_1\right)dx \nonumber\\
&= &\int_{\mathbb{R}}e^{yx_1}\omega(x_1)^{-q}\left(\int_{\mathbb{R}} e^{y(x-x_1)}\omega(x-x_1)^{-q} dx\right)dx_1 = \zeta(y)^2.
\end{eqnarray}
So, $0\leq \zeta(y)\leq1$.\\
Now define $\xi:[0,\rho_{2,\omega}q)\to(0,\infty)$ by $$\xi(y)=\int_0^\infty e^{yx}\omega(x)^{-q} \quad (y\in[0,\rho_{2,\omega}q)).$$ Then $0\leq\xi\leq\zeta\leq1$ and $\xi$ is an increasing continuous function. So, it is uniformly continuous on $[0,\rho_{2,\omega}q)$ and can be extended to $[0,\rho_{2,\omega}q]$ giving $$\xi(\rho_{2,\omega})=\int_0^\infty e^{\rho_{2,\omega}qx}\omega(x)^{-q}dx<\infty.$$ Now, \begin{align*} \int_\mathbb{R} e^{\rho_{2,\omega}qx}\omega(x)^{-q}dx &= \int_{-\infty}^0 e^{\rho_{2,\omega}qx}\omega(x)^{-q}dx+\int_0^\infty e^{\rho_{2,\omega}qx}\omega(x)^{-q}dx \\ &\leq \int_{-\infty}^0 \omega(x)^{-q}dx+\int_0^\infty e^{\rho_{2,\omega}qx}\omega(x)^{-q}dx \\ &<\infty. \end{align*} Thus, if $\rho_{2,\omega}>0$, then the function $h(x)=e^{\rho_{2,\omega}x}$ is in $L^q(\mathbb{R},\omega^{-1})$ and so for any $z=\rho_{2,\omega}+ib\in\mathbb{C}$, \begin{align*}\int_\mathbb{R}|f(x)||e^{zx}|dx & \leq \int_\mathbb{R}|f(x)|e^{\rho_{2,\omega}x} \omega(x) \omega(x)^{-1} dx \\ & \leq \left(\int_\mathbb{R}|f(x)|^p\omega(x)^p dx \right)^{1/p} \left(\int_\mathbb{R}e^{\rho_{2,\omega}qx} \omega(x)^{-q} dx \right)^{1/q}\\ & <\infty. \end{align*} Similarly, if $\rho_{1,\omega}<0$, then $g(x)=e^{\rho_{1,\omega}x}$ is in $L^q(\mathbb{R},\omega^{-1})$ and $\int_\mathbb{R}|f(x)||e^{zx}|dx <\infty$ for $z=\rho_{1,\omega}+ib\in\mathbb{C}$. So, $\int_\mathbb{R}f(x)e^{zx}dx $ is well defined for all $z\in S_\omega$.\\

For $z\in S_\omega$, let $\phi_z:F\to\mathbb{C}$ be defined as $$\phi_z(f)=\int_\mathbb{R}f(x)e^{zx}dx \quad (f\in F).$$ Then it is clear that $\phi_z\in\Delta(F)$. Let $\phi\in\Delta(F)$. Then using standard techniques we get a continuous function $\alpha:\mathbb{R}\to\mathbb{C}$ satisfying $\alpha(x+y)=\alpha(x)\alpha(y)$ and $|\alpha(x)|\leq\omega(x)$ for all $x,y\in\mathbb{R}$ such that $$\phi(f)=\int_\mathbb{R}f(x)\overline{\alpha(x)}dx \quad (f\in F).$$ Since $\alpha$ is continuous and $\alpha(x+y)=\alpha(x)\alpha(y)\ (x,y\in\mathbb{R})$, there is some $z=a+ib\in\mathbb{C}$ such that $\alpha(x)=e^{-izx}\ (x\in\mathbb{R})$. Since $|\alpha(x)|\leq\omega(x)$ for all $x\in\mathbb{R}$, we get $e^{bx}\leq\omega(x)$ for all $x\in\mathbb{R}$. Taking $z_0=b+ia$, we get $z_0\in S_\omega$ from last inequality and $\phi=\phi_{z_0}$ as for any $f\in F$, \begin{align*} \phi(f)=\int_\mathbb{R}f(x)\overline{e^{-izx}}dx =\int_\mathbb{R}f(x){e^{z_0x}}dx=\phi_{z_0}(f). \end{align*} Thus $\Delta(F)=\{\phi_z:z\in S_\omega\}$. Using calculation as in inequality (\ref{EQQ}) we get $\int_{\mathbb{R}}|e^{zx}|^q\omega(x)^{-q}dx\leq 1$ for $z \in S_\omega$.
\end{proof}

Let $f\in L^p(\mathbb{R},\omega,\mathcal{A})$. The fact that $\int_\mathbb{R}f(x)e^{zx}dx $ is well defined for all $z\in S_\omega$ can be shown similarly as shown in the last Proposition by replacing $|f(x)|$ by $\|f(x)\|$. The \emph{Fourier transform} $\widehat{f}:S_\omega \to \mathcal{A}$ of $f$ is defined as $$\widehat{f}(z)=\int_\mathbb{R} f(x) e^{zx}dx\quad(z \in S_\omega).$$

\begin{lem}\label{lmrp1}
Let $\omega$ be a weight on $\mathbb{R}$.
\begin{enumerate}
\item If $\rho_{2,\omega}>0$ and $\gamma \in (0,1)$, then there is a positive constant $K$ such that $e^{x^\gamma}\leq K\omega(x)$ for all $x\geq0$.
\item If $\rho_{1,\omega}<0$ and $\gamma \in (0,1)$, then there is a positive constant $K$ such that $e^{-x^\gamma}\leq K\omega(x)$ for all $x\leq0$.
\item If $\rho_{1,\omega}<0<\rho_{2,\omega}$ and $\gamma \in (0,1)$, then there is a positive constant $K$ such that $e^{|x|^\gamma}\leq K\omega(x)$ for all $x \in \mathbb{R}$.
\end{enumerate}
\end{lem}
\begin{proof}
Take $\gamma\in(0,1)$. Since $\rho_{2,\omega}>0$ and $\lim_{x\to\infty} \log (e^{x^\gamma})^{1/x}=0$, there is $R>0$ such that $\log (e^{x^\gamma})^{1/x} \leq \frac{\rho_{2,\omega}}{2}\leq \log \omega(x)^{1/x}$ for all $x>R$. This implies that $e^{x^\gamma}\leq\omega(x)$ for all $x>R$. Since $[0,R]$ is compact in $\mathbb{R}$, there is some $K\geq1$ such that $e^{x^\gamma}\leq K\omega(x)$ for all $x\in[0,R]$. Therefore $e^{x^\gamma}\leq K\omega(x)$ for all $x>0$. The proof of (ii) is similar and (iii) follows from (i) and (ii).
\end{proof}

\begin{lem} \label{lmrp2}
Let $\omega$ be an almost monotone algebra weight on $\mathbb{R}$.
\begin{enumerate}
\item If $\rho_{2,\omega}=0$, define $\widetilde{\omega}:\mathbb{R} \to [1,\infty)$ by $\widetilde{\omega}(x)=\sup \{\omega(h):0\leq h \leq x \}$ if $x\geq0$ and $\widetilde{\omega}(x)=\omega(x)$ if $x<0$.
\item If $\rho_{1,\omega}=0$, define $\widetilde{\omega}:\mathbb{R} \to [1,\infty)$ by $\widetilde{\omega}(x)=\sup \{\omega(h):x\leq h \leq 0 \}$ if $x\leq0$ and $\widetilde{\omega}(x)=\omega(x)$ if $x>0$.
\end{enumerate}
Then in any case $\widetilde{\omega}$ is an almost monotone algebra weight on $\mathbb{R}$ and there is a positive constant $K$ such that $\widetilde{\omega}(x)\leq K \omega(x)$ for all $x \in \mathbb{R}$.
\end{lem}
\begin{proof}
Since, by \cite[Lemma 1.3.3]{Kan}, $\omega$ is bounded on compact subsets of $\mathbb{R}$, $\widetilde{\omega}$ is well-defined. Let $x,y\in\mathbb{R}$. If both $x$ and $y$ are non-negative (or negative), then $\widetilde\omega(x+y)\leq\widetilde\omega(x)\widetilde\omega(y)$. Let $x\geq0$ and $y<0$. If $x+y\geq0$, then $\widetilde\omega(x+y) \leq \widetilde\omega(x) \leq \widetilde\omega(x)\omega(y) = \widetilde\omega(x)\widetilde\omega(y)$ as $\widetilde\omega$ is increasing on $[0,\infty)$ and $y<0$. If $x+y<0$, then $\widetilde\omega(x+y) = \omega(x+y) \leq \omega(x)\omega(y) \leq \widetilde\omega(x)\widetilde\omega(y)$. Therefore $\widetilde\omega$ is a weight. Since $\rho_{2,\omega}=0$ and $\omega$ is an almost monotone algebra weight, there is $K>0$ such that $\omega(x)\leq K\omega(x+(y-x))\leq K\omega(y)$ for all $0\leq x\leq y$. Let $x\geq0$. Then $\omega(h)\leq K\omega(x)$ for all $0\leq h\leq x$. This implies that $\widetilde\omega(x)\leq K\omega(x)$. It is clear that $K\geq1$. So, $\widetilde \omega(x)\leq K\widetilde\omega(x)=K\omega(x)$ for all $x<0$. Since $\omega(x)\leq \widetilde\omega(x)\leq K\omega(x)$ for all $x\geq0$ and $\rho_{2,\omega}=0,\ \rho_{2,\widetilde\omega}=0$. The inequality $\omega(x)\leq \widetilde\omega(x)\leq K\omega(x)$ for all $x\in\mathbb{R}$ implies that $\int_{\mathbb{R}} \widetilde\omega^{-q}<\infty$. Thus $\widetilde\omega$ is an almost monotone algebra weight. The proof of (ii) is similar.
\end{proof}

\begin{lem}\label{lmcrp}
Let $1<p<\infty, q$ be its conjugate index, $\omega$ be an almost monotone algebra weight on $\mathbb{R}, \mathcal{A}$ be a unital Banach algebra, $F=\{\widehat{f}:f\in L^p(\mathbb{R},\omega,\mathbb{C})\}$, and let $\mathcal{F}=\{\widehat{f}:f\in L^p(\mathbb{R},\omega,\mathcal{A})\}$. Then $F$ and $\mathcal{F}$ satisfies all conditions from (i) to (v) in \textbf{($\bigstar$)}.
\end{lem}
\begin{proof}
Let $\phi\in\Delta(F)$, and let $x=a_1f_1+a_2f_2+\dots+a_nf_n\in\mathcal{F}$ for some $n\in\mathbb{N}$ where $a_i\in\mathcal{A}$ and $f_i\in F$ for $1\leq i\leq n$. Since the Gel'fand space $\Delta(F)$ of $F$ is identified with the strip $S_\omega$ via the map $\phi\mapsto\phi_z$ where $\phi_z(f)=\widehat f(z) \ (f\in F)$, there is some $z=a+ib\in S_\omega$ such that $\phi=\phi_z$. So, \begin{eqnarray*} &&\|a_1\phi_z(f_1)+a_2\phi_z(f_2)+\dots+a_n\phi_z(f_n)\| \\ &=& \|a_1\widehat{f_1}(z)+a_2\widehat{f_2}(z)+\dots+a_n\widehat{f_n}(z)\| \\ &\leq& \int_{\mathbb{R}} \|a_1f_1(x)+a_2f_2(x)+\dots+a_nf_n(x)\||e^{zx}|dx \\ &\leq& \int_{\mathbb{R}} \|a_1f_1(x)+a_2f_2(x)+\dots+a_nf_n(x)\|\omega(x)e^{ax}\omega(x)^{-1}dx \\ &\leq& \left(\int_{\mathbb{R}} \|a_1f_1(x)+a_2f_2(x)+\dots+a_nf_n(x)\|^p\omega(x)^p dx\right)^{1/p} \left(\int_\mathbb{R}e^{aqx}\omega(x)^{-q}dx \right)^{1/q} \\ &\leq& \|x\|_{\mathcal{F}}. \end{eqnarray*}
\end{proof}

\begin{proof}[\textbf{Proof of Theorem \ref{thrp}}]
As done in Theorem \ref{thr} we get $M>0$ and $r_1,r_2\in\mathbb{R}$ such that $\rho_{1,\omega}\leq r_1 \leq 0 \leq r_2 \leq \rho_{2,\omega}$ and $\mathbf  e+\widehat{f}(z)$ is left invertible in $\mathcal{A}$ for all $z\in\{x+iy\in S_\omega:|y|>M\}\cup\{x+iy:r_1\leq x\leq r_2,-M\leq y\leq M\}$. Notice that $r_1=0$ if $\rho_{1,\omega}=0$, $\rho_{1,\omega}<r_1<0$ if $\rho_{1,\omega}<0$, $r_2=0$ if $\rho_{2,\omega}=0$ and $0<r_2<\rho_{2,\omega}$ if $\rho_{2,\omega}>0$. Define $\nu:\mathbb{R}\to [1,\infty)$ as follows:\\
If $\rho_{1,\omega}=0=\rho_{2,\omega}$, then take $\nu=\omega$.\\
Let $\rho_{1,\omega}=0$ and $\rho_{2,\omega}>0$. Then $\mathbf  e+\widehat f(x+iy)$ is left invertible in $\mathcal{A}$ for all $x\in [0,r_2]$ and $y\in [-M,M]$. By lemma \ref{lmrp2}, there is some $K_1>0$ such that $\omega(x)\leq \widetilde\omega(x)\leq K_1\omega(x)$ for all $x\leq0$, where $\widetilde\omega(x)=\sup\{\omega(h):x\leq h \leq0\}$ for all $x\leq0$. Take $\gamma\in(0,1)$. Then, by lemma \ref{lmrp1}, there is $K_2>0$ such that $e^{x^\gamma}\leq K_2\omega(x)$ for all $x>0$. Let $K=\max\{K_1,K_2\}$ and define $$\nu(x)=\begin{cases} \widetilde\omega(x) & x<0 \\ e^{\frac{r_2x+x^\gamma}{2}} & x\geq0 \end{cases}.$$
Let $\rho_{1,\omega}<0$ and $\rho_{2,\omega}=0$. If $r_1\leq x\leq0$ and $-M\leq y\leq M$, then $\mathbf  e+\widehat f(z)$ is left invertible in $\mathcal{A}$ for all $z=x+iy$. Take $\gamma\in(0,1)$. Then, by lemma \ref{lmrp1}, there is $K_1>0$ such that $e^{-x^\gamma}\leq K_1\omega(x)$ for all $x<0$. By lemma \ref{lmrp2}, there is is some $K_2>0$ such that $\omega(x)\leq \widetilde\omega(x)\leq K_2\omega(x)$ for all $x\geq0$, where $\widetilde\omega(x)=\sup\{\omega(h):0\leq h\leq x\}$ for all $x\geq0$. Take $K=\max\{K_1,K_2\}$ and define $$\nu(x)=\begin{cases} e^{\frac{r_1x+x^\gamma}{2}} & x<0 \\ \widetilde\omega(x) & x\geq0 \end{cases}.$$
Let $\rho_{1,\omega}<0$ and $\rho_{2,\omega}>0$. Then $\mathbf  e+\widehat f(z)$ is left invertible in $\mathcal{A}$ for all $z=x+iy$ with $r_1\leq x \leq r_2$ and $-M\leq y\leq M$. Take $\gamma\in(0,1)$. Then, by lemma \ref{lmrp1}, there is $K>0$ such that $e^{|x|^\gamma}\leq K\omega(x)$ for all $x\in\mathbb{R}$. Define $$\nu(x)=\begin{cases} e^{\frac{r_1x+x^\gamma}{2}}  & x<0 \\ e^{\frac{r_2x+x^\gamma}{2}}  & x\geq0 \end{cases}.$$
Then $\nu$ is an almost montone algebra weight on $\mathbb{R}$ such that $1\leq\nu\leq K\omega$. This gives $L^p(\mathbb{R},\omega,\mathcal{A})\subset L^p(\mathbb{R},\nu,\mathcal{A})$ and so $f\in L^p(\mathbb{R},\nu,\mathcal{A})$. Let $F=\{\widehat{g}:g\in L^p(\mathbb{R},\nu,\mathbb{C})\}$, and let $\mathcal{F}=\{\widehat{g}:g\in L^p(\mathbb{R},\nu,\mathcal{A})\}$. Then $\widehat f\in \mathcal{F}$. The strip $S_\nu=\{x+iy:\sqrt r_1\leq x \leq \sqrt r_2,\ y\in \mathbb{R}\}$ is the Gel'fand space $\Delta(F)$ of $F$ by the map $z\mapsto\phi_z$, where $\phi_z(g)=g(z) \ (g\in F).$ Every generated homomorphism of $\mathcal{F}$ is of the form $\Phi_z$ for $z\in S_\nu$, where $\Phi_z(g)=g(z) \ (g\in\mathcal{F})$. Since $\mathbf  e+\Phi_z(\widehat f)=\mathbf  e+\widehat f(z)$ is left invertible in $\mathcal{A}$, by Lemma \ref{lmcrp} and Theorem \ref{th1}, $\mathbf  e+\widehat{f}$ has a left inverse in $\mathcal{F}$, say, $\mathbf  e+\widehat{g}$. This implies that $g\in L^p(\mathbb{R},\nu,\mathcal{A})$ is such that $\mathbf{1}+g$ is a left inverse of $\mathbf{1}+f$. This completes the proof.
\end{proof}

\section{Weighted off diagonal decay of infinite matrices of operators}
Let $X$ be a complex Banach space and $B(X)$ be the Banach algebra of all bounded linear operators on $X$. A map $R: \mathbb{Z} \to B(X)$ is a \emph{resolution of the identity} if $R(n)$ is a projection on $X$ for all $n$, $R(k)R(l)=0$ for all $k \neq l$ and the series $\sum_{k \in \mathbb{Z}} R(k)x$ converges to $x$ unconditionally for all $x \in X$.
Following \cite{Ba}, we make following assumptions:
\begin{enumerate}
\item $$1 \leq C_R=\sup_{\alpha_k \in \mathbb{T}} \| \sum_{k \in \mathbb{Z}} \alpha_k R(k) \| < \infty.$$
\item There is a positive constant $M_R$ such that $$\left\| \sum_{1 \leq i \leq n} R(i+k)AR(i) \right\| \leq M_R \max_{1 \leq i \leq n} \| R(i+k)AR(i) \|\;(k \in \mathbb{Z}, n \in \mathbb{N}, A \in B(X)).$$
\end{enumerate}

For each $A \in B(X)$, consider the matrix $(A_{ij})_{i,j \in \mathbb{Z}}$, where $A_{ij}=R(i)AR(j) \in B(X)$. Then for each $k \in \mathbb{Z}$, the \emph{$k^{th}$- diagonal} of the matrix of $A$ is the collection $\{A_{ij}: i-j=k\} \subset B(X)$. Define $d_A: \mathbb{Z} \to [0,\infty)$ by $$d_A(k)=\sup_{i-j=k}\|A_{ij}\| \ \ (k \in \mathbb{Z}).$$ Then $d_A$ characterizes the off-diagonal decay of the entries of $(A_{ij})$.

Consider the representation $T:\mathbb{T} \to B(X)$ defined by $$T(e^{it})x=\sum_{n \in \mathbb{Z}}R(n)xe^{int} \quad (x \in X, e^{it} \in \mathbb{T}).$$ Then from Assumption (i), $\|T(e^{it})\| \leq C_R$ for all $e^{it} \in \mathbb{T}$ and so $T$ is a bounded strongly continuous representation.

Now, for $A \in B(X)$, define \cite{Ba} $\widehat{A}: \mathbb{T} \to B(X)$ as $$\widehat{A}(e^{it})=T(e^{it})AT(e^{-it}) \quad (e^{it} \in \mathbb{T}).$$ Let $$\widehat{A}(e^{it})(x) \thicksim \sum_{n \in \mathbb{Z}}A_nxe^{int} \quad (x \in X),$$ be the Fourier series of $\widehat{A}$, where $A_n \ (n \in \mathbb{Z})$ are the Fourier coefficients given by $$A_n=\frac{1}{2\pi} \int_0^{2\pi} \widehat{A}(e^{it})e^{-int} dt.$$ The series $\sum_{n \in \mathbb{Z}} A_n$ is the \emph{Fourier series} of $A$ and $A_n \in B(X) \ (n \in \mathbb{Z})$ are the \emph{Fourier coefficients} of $A$. By \cite[Lemma 1]{Ba} $$A_n=\sum_{i-j=n}R(i)AR(j) =\sum_{i-j=n} A_{ij} \quad(n \in \mathbb{Z}),$$ there are positive constants $C_1(R)$ and $C_2(R)$ such that $$C_1(R)d_A(n) \leq \|A_n\| \leq C_2(R)d_A(n) \quad(n \in \mathbb{Z}),$$ and by \cite[Corollary to Lemma 1]{Ba}, all diagonals entries of $A_n$ are $0$ except the $n^{th}$- diagonal which is same as the $n^{th}$- diagonal of $A$ (Assumption (ii) is used in it). By \cite[Lemma 2]{Ba}, $$B_\omega(X)=\{ A \in B(X) : \sum_{n \in \mathbb{Z}} d_A(n) \omega(n) < \infty \}$$ is a subalgebra of $B(X)$ and by \cite[Corollary to Lemma 2]{Ba}, it is a Banach algebra with the norm $$\|A\|_\omega = \sum_{n \in \mathbb{Z}} d_A(n) \omega(n) \quad (A \in B_\omega(X)).$$ In \cite{Ba}, Baskakov proved that if $\omega$ is an admissible weight and $A\in B_\omega(X)$ is invertible in $B(X)$, then $A^{-1}\in B_\omega(X)$.\\
 Let $0<p<\infty$, and let $\omega$ be a weight on $\mathbb Z$. Let $$B_{p\omega}(X)=\{A \in B(X) :\sum_{n \in \mathbb{Z}} d_A(n)^p \omega(n)^p < \infty \}$$

\begin{lem} \label{lm3}
Let $0<p<\infty$, $\omega$ be a weight on $\mathbb Z$, and let $X$ be a Banach space.
\begin{enumerate}
\item If $0<p\leq 1$, then $B_{p\omega}(X)$ is a $p$-Banach algebra with the $p$-norm $$\|A\|_{p\omega}=\sum_{n\in \mathbb Z}d_A(n)^p\omega(n)^p\quad(A \in B_{p\omega}(X)).$$
\item If $1<p<\infty$ and $\omega$ satisfies $\omega^{-q}\star \omega^{-q}\leq \omega^{-q}$, where $q$ is the conjugate index of $p$, then $B_{p\omega}(X)$ is a Banach algebra with the norm $$\|A\|_{p\omega}=\left(\sum_{n\in \mathbb Z}d_A(n)^p\omega(n)^p\right)^{\frac{1}{p}}\quad(A \in B_{p\omega}(X)).$$
\end{enumerate}
\end{lem}

\begin{thm} \label{tho1}
Let $0<p\leq 1$. Let $X$ be a complex Banach space, $\omega$ be a weight on $\mathbb{Z}$, and let  $A\in B_{p\omega}(X)$ be invertible in $B(X)$. Then there is a weight $\nu$ on $\mathbb{Z}$ such that $1 \leq \nu \leq \omega$, $\nu$ is constant if and only if $\omega$ is constant and $B=A^{-1} \in B_{p\nu}(X)$. In particular, if $\omega$ is an admissible weight, then $A^{-1} \in B_{p\omega}(X)$.
\end{thm}
\begin{proof}
Consider the function $f:\mathbb{Z} \to B(X)$ defined by $f(n)=A_n\;(n \in \mathbb Z)$, where $A_n$ are Fourier coefficients of $A$. Then by hypothesis $f \in \ell^p(\mathbb{Z},\omega,B(X))$ as $$\sum_{n\in \mathbb{Z}} \|f(n)\|^p \omega(n)^p = \sum_{n\in \mathbb{Z}} \|A_n\|^p \omega(n)^p \leq C_2(R)^p \sum_{n\in \mathbb{Z}} d_A(n)^p \omega(n)^p < \infty.$$ Let $\mathcal{A}=B(X)$. By \cite[Theorem 2]{Ba}, the Fourier transform $\widehat{f}: \mathbb{T} \to B(X)$ of $f$ is given by $$\widehat{f}(e^{it})= \widehat{A}(e^{-it}) \quad (e^{it} \in \mathbb{T})$$ and $\widehat{f} \in \mathcal{A}_p(\omega)$, as defined in the proof of Theorem \ref{thm1}. Since $\widehat{f}(e^{it}) = \widehat{A}(e^{-it}) = T(e^{-it})AT(e^{it})$ is invertible in $B(X)$ for all $e^{it} \in \mathbb{T}$, by Corollary \ref{crm1}, there is a weight $\nu$ on $\mathbb{Z}$ such that $1 \leq \nu \leq \omega$, $\nu$ is constant if and only if $\omega$ is constant and $f^{-1} \in \ell^p(\mathbb{Z},\nu,B(X))$. Now, $$\widehat{f^{-1}}(e^{it}) = \widehat{A}(e^{-it})^{-1} = T(e^{-it})BT(e^{it}) \quad(e^{it} \in \mathbb{T}).$$ So, $$\widehat{B} (e^{it}) = T(e^{it})BT(e^{-it}) = \sum_{n \in \mathbb{Z}} B_n e^{int} = \sum_{n \in \mathbb{Z}} f^{-1}(n) e^{-int} \quad (e^{it}\in\mathbb{T}).$$ Thus $B=A^{-1} \in B_{p\nu}(X)$.
\end{proof}

\begin{thm} \label{tho2}
Let $1<p<\infty$, $q$ be the conjugate index of $p$, and let $\omega$ be an almost monotone algebra weight on $\mathbb Z$. Let $X$ be a complex Banach space, and let  $A\in B_{p\omega}(X)$ be invertible in $B(X)$. Then there is an almost monotone algebra weight $\nu$ on $\mathbb{Z}$ such that $1 \leq \nu \leq K\omega$ for some $K>0$ and $B=A^{-1} \in B_{p\nu}(X)$. In particular, if $\omega$ is an admissible weight, then $A^{-1} \in B_{p\omega}(X)$.
\end{thm}
\begin{proof}
The proof is similar to that of  Theorem \ref{tho1} by using Corollary \ref{crm2} in place of Corollary \ref{crm1}.
\end{proof}

\begin{remark}
We do not know whether there is a weight $\omega$ on $\mathbb Z$ satisfying $\omega^{-q}\star \omega^{-q}\leq \omega^{q}$, $\sum_{n\in \mathbb Z}\omega^{-q}(n)<\infty$ for some $1<q<\infty$ and $\rho_{2\omega}=1$ but there is no $K>0$ satisfying $\omega(n)\leq K\omega(n+k)$ for all $n,k \in \mathbb N_0$.
\end{remark}

\subsection*{Acknowledgements}
The first author would like to thank SERB, India, for the MATRICS grant no. MTR/2019/000162. The second author gratefully acknowledges Junior Research Fellowship (NET) from CSIR, India.


\end{document}